\documentclass{article}

\title{Tensor-generating line bundles on toric varieties}
\renewcommand\footnotemark{}
\thanks{2020 {\em Mathematics Subject Classification.} 14M25 (primary), 14F08, 52B20 (secondary)}\thanks{{\em Key words and phrases.} toric variety, derived category, tensor generation, big line bundle}
\author{Daigo Ito, Michael R. Zeng, and Xiangru Zeng}
\date{}

\makeatletter
\IfFileExists{paper/style.tex}
  {}
  {\def\input@path{{../}}}
\makeatother

\usepackage{stix} 
\usepackage[english]{babel} 
\usepackage{amsmath,amsfonts,amsthm,amscd}
\usepackage{ascmac}
\usepackage{appendix}
\usepackage{bm}
\usepackage{cases}
\usepackage{changepage}
\usepackage{enumitem}
\usepackage{etoolbox}
\usepackage{float}
\usepackage{geometry}
\usepackage{graphicx}
\usepackage[utf8]{inputenc}
\usepackage [autostyle, english = american]{csquotes}
\usepackage{scalerel}
\usepackage{ytableau}
\usepackage{tikz-cd} 
\usepackage{tocloft}
\usepackage{setspace}
\usepackage[color,all]{xy}
\usepackage[cal=euler]{mathalfa}
\usepackage{url}
\usepackage{CJKutf8}
\makeatletter
\newcommand{\@endstuff}{\par\vspace{\baselineskip}\noindent\small
\begin{tabular}{@{}l}\scshape{Daigo Ito} \\ \scshape Department of Mathematics, Columbia University, 2990 Broadway, New York, NY 10021\\\textrm{E-mail address:} \url{di2260@columbia.edu}\\\textrm{Webpage:} \url{https://daigoi.github.io} \end{tabular}
\par\vspace{\baselineskip}\noindent\small
\begin{tabular}{@{}l}\scshape{Michael Ruofan Zeng} \\ \scshape Department of Mathematics, University of Washington, 4110 E Stevens Way NE, Seattle, WA 98105 \\\textrm{E-mail address:} \url{zengrf@uw.edu} \\ \textrm{Webpage:} \url{https://zengrf.github.io} \end{tabular}
\par\vspace{\baselineskip}\noindent\small
\begin{tabular}{@{}l}\scshape{Xiangru Zeng} \\ \scshape Department of Mathematics, University of California, Berkeley, Evans Hall, CA 94720-3840 \\\textrm{E-mail address:} \url{xiangru_zeng@berkeley.edu} \end{tabular}}
\AtEndDocument{\@endstuff}
\makeatother

\usepackage{hyperref}
\usepackage[alphabetic,backrefs]{amsrefs}
\hypersetup{colorlinks,linkcolor=blue,citecolor=red}

\usepackage{cleveref}

\DeclareMathOperator{\aut}{Aut}

\let \hom \relax
\DeclareMathOperator{\hom}{Hom}

\let \im \relax
\DeclareMathOperator{\im}{Im}

\DeclareMathOperator{\perf}{\mathsf{Perf}}

\DeclareMathOperator{\pic}{Pic}

\DeclareMathOperator{\spec}{Spec}

\DeclareMathOperator{\supp}{Supp}

\newcommand {\bb}{\mathbb}
\renewcommand {\cal}{\mathcal}
\newcommand{\ecal}{\mathscr}
\renewcommand {\epsilon}{\varepsilon}

\newcommand {\bra}[1]{\langle{#1}\rangle}

\newcommand {\emp}{\emptyset}

\newcommand {\inv}{^{-1}}
\newcommand {\iso}{\cong}

\newcommand {\tens}{\otimes}

\newsavebox{\pullbacks}
\sbox\pullbacks{%
\begin{tikzpicture}%
\draw (0,0) -- (1ex,0ex);%
\draw (1ex,0ex) -- (1ex,1ex);%
\end{tikzpicture}}

\newcommand{\D}{\mathsf{D}}

\newcommand{\h}{{\mathrm{H}}}

\newcommand{\injar}{\ar@{^(->}} 
\newcommand{\prarrow}[2]{\ar@<0.5ex>[r]^-{#1} \ar@<-0.5ex>[r]_-{#2}}
\newcommand{\plarrow}[2]{\ar@<0.5ex>[l]^-{#1} \ar@<-0.5ex>[l]_-{#2}}
\newcommand{\pdarrow}[2]{\ar@<0.5ex>[d]^-{#1} \ar@<-0.5ex>[d]_-{#2}}
\newcommand{\puarrow}[2]{\ar@<0.5ex>[u]^-{#1} \ar@<-0.5ex>[u]_-{#2}}

\theoremstyle{plain}
\newtheorem{theorems}{Theorem}[section] 
\newtheorem{claims}{Claim}[theorems]
\newtheorem{conjectures}[theorems]{Conjecture}
\newtheorem{corollaries}[theorems]{Corollary}
\newtheorem{lemmas}[theorems]{Lemma} 
\newtheorem{props}[theorems]{Proposition}

\newtheorem{penmdef}[claims]{Definition}

\theoremstyle{definition}
\newtheorem{constructions}[theorems]{Construction}
\newtheorem{definitions}[theorems]{Definition} 
\newtheorem{axioms}[theorems]{Axiom} 
\newtheorem{examples}[theorems]{Example}     
\newtheorem{notations}[theorems]{Notation}         

\newtheorem{question}[theorems]{Question}
\newtheorem{obss}[theorems]{Observation}
\newtheorem{penmlem}[claims]{Lemma}
\newtheorem{penmthm}[claims]{Theorem}
\newtheorem{penmcor}[claims]{Corollary}
\newtheorem{penmeg}[claims]{Example} 
\newtheorem{inclaims}[claims]{Claim}

\theoremstyle{remark}
\newtheorem{remarks}[theorems]{Remark}

\newtheorem{penmrem}[claims]{Remark}

\makeatletter
\newtheoremstyle{indented}
  {1pt}
  {1pt}
  {\addtolength{\@totalleftmargin}{1.5em}
   \addtolength{\linewidth}{-1.5em}
   \parshape 1 1.5em \linewidth}
  {}
  {\bfseries}
  {.}
  {.5em}
  {}
\makeatother

\theoremstyle{indented}
\newtheorem{pinddef}[claims]{Definition} 
\newtheorem{pindlem}[claims]{Lemma}
\newtheorem{pindthm}[claims]{Theorem}
\newtheorem{pindcor}[claims]{Corollary}
\newtheorem{pindeg}[claims]{Example} 
\newtheorem{pindrem}[claims]{Remark} 
\newtheorem{pindq}[claims]{Question} 
\newtheorem{pindc}[claims]{Conjecture} 
\newtheorem{pindclaim}[claims]{Claim}

\newenvironment{theorem}
{
	\pushQED{\qed}\begin{theorems}}
	{\popQED\end{theorems}}

\newenvironment{prop}
{
	\pushQED{\qed}\begin{props}}
	{\popQED\end{props}}

\newenvironment{notation}
{
	\pushQED{\qed}\begin{notations}}
	{\popQED\end{notations}}

\newenvironment{corollary}
{
	\pushQED{\qed}\begin{corollaries}}
	{\popQED\end{corollaries}}

\newenvironment{definition}
{
	\pushQED{\qed}\begin{definitions}}
	{\popQED\end{definitions}}

\newenvironment{lemma}
{
	\pushQED{\qed}\begin{lemmas}}
	{\popQED\end{lemmas}}

\newenvironment{remark}
{
	\pushQED{\qed}\begin{remarks}}
	{\popQED\end{remarks}}
	
\newenvironment{example}
{
	\pushQED{\qed}\begin{examples}}
	{\popQED\end{examples}}

\newenvironment{construction}
{
	\pushQED{\qed}\begin{constructions}}
	{\popQED\end{constructions}}

\theoremstyle{plain}

\theoremstyle{definition}

\crefformat{equation}{(#2#1#3)}
\Crefformat{equation}{(#2#1#3)}
\crefname{thmx}{Theorem}{Theorems}          \Crefname{thmx}{Theorem}{Theorems}
\crefname{defx}{Definition}{Definitions}    \Crefname{defx}{Definition}{Definitions}
\crefname{theorems}{Theorem}{Theorems}      \Crefname{theorems}{Theorem}{Theorems}
\crefname{props}{Proposition}{Propositions} \Crefname{props}{Proposition}{Propositions}
\crefname{lemmas}{Lemma}{Lemmas}            \Crefname{lemmas}{Lemma}{Lemmas}
\crefname{corollaries}{Corollary}{Corollaries} \Crefname{corollaries}{Corollary}{Corollaries}
\crefname{definitions}{Definition}{Definitions} \Crefname{definitions}{Definition}{Definitions}

\crefname{constructions}{Construction}{Constructions} \Crefname{constructions}{Construction}{Constructions}
\crefname{examples}{Example}{Examples}      \Crefname{examples}{Example}{Examples}
\crefname{remarks}{Remark}{Remarks}         \Crefname{remarks}{Remark}{Remarks}
\crefname{question}{Question}{Questions}    \Crefname{question}{Question}{Questions}
\crefname{notations}{Notation}{Notations}   \Crefname{notations}{Notation}{Notations}
\crefname{conjectures}{Conjecture}{Conjectures} \Crefname{conjectures}{Conjecture}{Conjectures}
\crefname{appendix}{Appendix}{Appendices}   \Crefname{appendix}{Appendix}{Appendices}

\MakeOuterQuote{"}

\DeclareTextAccent{\u}{T1}{8}

\usetikzlibrary{arrows.meta,calc}
\colorlet{clsgray1}{black!14}
\colorlet{clsgray2}{black!26}
\colorlet{clsgray3}{black!39}
\colorlet{clsgray4}{black!52}
\colorlet{clsface1}{black!6}
\colorlet{clsface2}{black!14}
\tikzset{
    cls plot/.style={
        line cap=round, line join=round, text=black,
        ray/.style={-{Latex[length=1.5mm,width=1mm]}, draw=black, line width=0.55pt},
        tested ray/.style={line width=0.85pt},
        gen/.style={circle, fill=black, inner sep=0.95pt},
        origin/.style={circle, fill=black!60, inner sep=0.8pt},
        lattice/.style={circle, fill=black!22, inner sep=0.6pt},
        lbl/.style={font=\footnotesize, text=black, inner sep=2pt},
        slbl/.style={lbl},
        title/.style={lbl},
        slice/.style={dash pattern=on 2pt off 1.6pt, draw=black!65, line width=0.4pt},
        chord/.style={slice, draw=black!85, line width=0.45pt},
        sep/.style={draw=black, line width=0.45pt},
        edgeBase/.style={draw=black!55, line width=0.4pt},
        edgeHidden/.style={edgeBase, dash pattern=on 2pt off 1.8pt},
        edgeBlow/.style={draw=black, line width=0.95pt},
        edgeBlowHidden/.style={edgeBlow, dash pattern=on 2pt off 1.8pt},
        vtx/.style={gen},
        elbl/.style={lbl, inner sep=1pt}
    }
}
\newcommand{\CLSray}[2][]{%
    \draw[ray,#1] (O) -- ($(#2)!-2.5mm!(O)$);
}
\newcommand{\CLSwall}[3]{%
    \fill[#3] (O) -- (#1) -- (#2) -- cycle;
}

\begin{document}

\maketitle

\begin{abstract}
A line bundle is $\otimes$-generating if its powers generate the derived category of perfect complexes. Refining the derived Nakai--Moishezon criterion of Ito--Olander, we characterize $\otimes$-generating line bundles on proper toric varieties by a combinatorial criterion on the fan, leading to the notion of $\otimes$-generating fans characterizing $\otimes$-generating canonical bundles. Furthermore, we prove an analogue of Kleiman's criterion of ampleness for $\otimes$-generation.
As applications, we provide several examples including: (i) smooth proper fourfolds with no $\otimes$-generating line bundle, by showing that products of two or more positive-dimensional proper toric varieties admit a $\otimes$-generating line bundle only when they are projective, (ii) families of smooth non-projective proper threefolds with $\otimes$-generating canonical bundles. 
Finally, we study subvarieties obstructing $\otimes$-generation of canonical bundles and implications to autoequivalences of derived categories. 
\end{abstract}

\tableofcontents

\section{Introduction}\label{sec:introduction}

A line bundle $\ecal L$ on a noetherian scheme $X$ is called \emph{$\otimes$-generating} if the tensor-powers $\{\ecal L^{\tens m}\}_{m \in \bb Z}$ generate the derived category $\perf(X)$ of perfect complexes. This condition was first formally studied in \cite{ito2025polarizations}, as a derived category analogue of ampleness, to generalize the reconstruction theorem of Bondal and Orlov \cite{bondal_orlov_2001}, which claims that (i) a smooth projective variety $X$ with (anti-)ample canonical bundle is determined by its derived category, and 
(ii) every exact autoequivalence of $\perf(X)$ is standard, that is, generated by automorphisms, line bundle twists, and shifts. 
Later in \cite{ito2025derived}, Ito and Olander characterize $\otimes$-generation on proper $k$-schemes by a derived category analogue of the Nakai--Moishezon criterion: $\ecal L$ is $\otimes$-generating if and only if $\ecal L|_Z$ is big or anti-big for every integral closed subscheme $Z \subseteq X$ \cite{ito2025derived}*{Theorem 1.1}. In particular, ample and anti-ample line bundles are $\otimes$-generating. On a curve, this condition reduces to having nonzero degree. On a smooth projective surface, a big line bundle is $\otimes$-generating if and only if it has nonzero degree on every curve of negative self-intersection \cite{ito2025derived}*{Lemma 5.3}.

In this paper, we study this positivity property on toric varieties. We give a combinatorial criterion, a numerical characterization by curves, existence and non-existence results, and a reflexivity-type characterization of $\otimes$-generating canonical bundles with families of examples. This includes the first smooth proper varieties known to us with no $\otimes$-generating line bundle at all and non-projective smooth proper threefolds whose canonical bundle is $\otimes$-generating.

Throughout the introduction, $X = X_\Sigma$ is a proper toric variety over a field $k$ with cocharacter lattice $N$, character lattice $M$, primitive ray generators $u_\rho$ for $\rho \in \Sigma(1)$, invariant prime divisors $D_\rho$, and orbit closures $V(\tau)$ for $\tau \in \Sigma$. For an invariant $\bb Q$-Cartier divisor $D = \sum_\rho a_\rho D_\rho$ and a cone $\tau \in \Sigma$, the \emph{restriction polytope} is
\[
P_D(\tau) := \{m \in M_{\bb R} \mid \bra{m, u_\rho} = -a_\rho \ \forall \rho \in \tau(1), \ \bra{m, u_\rho} \geq -a_\rho \ \forall \rho \in \sigma(1), \ \sigma \succ \tau\}.
\]

\begin{theorem}[\Cref{thm:orbit-closure,thm:implementable}, \Cref{prop:hyperplane-test}]\label{thm:A}
    Let $X = X_\Sigma$ be a proper toric variety over a field $k$ and $D$ an invariant $\bb Q$-Cartier divisor with $lD$ Cartier, where $l$ is a positive integer. Then, the following are equivalent.
    \begin{enumerate}
        \item The line bundle $\ecal O_X(lD)$ is $\otimes$-generating.
        \item The restriction of $\ecal O_X(lD)$ to every orbit closure $V(\tau)$ is big or anti-big.
        \item For every cone $\tau \in \Sigma$, either $\dim P_D(\tau)$ or $\dim P_{-D}(\tau)$ equals $\dim X - \dim \tau$.
    \end{enumerate}
    If moreover all $a_\rho > 0$, these are equivalent to a separation condition on the rescaled generators $u_\rho / a_\rho$: 
    
    \begin{enumerate}
        \item[(iv)] for every cone $\tau \neq \{0\}$, some affine hyperplane contains the points $u_\rho/a_\rho$ with $\rho \in \tau(1)$ and has the points $u_\rho/a_\rho$, for the remaining rays of the cones containing $\tau$, strictly on one side. \qedhere 
    \end{enumerate}
    
\end{theorem}

The equivalence of (i) and (ii) reduces the derived Nakai--Moishezon criterion of \cite{ito2025derived} to the finitely many invariant subvarieties. Conditions (iii) and (iv) are combinatorial and checkable as linear programs, for which we provide a SageMath implementation (see \hyperref[par:code-availability]{Code availability} below). We then specialize the combinatorial criterion to canonical bundles. For the anticanonical divisor $-K_X = \sum_\rho D_\rho$, the criterion suggests the following new class of complete fans which generalizes reflexive fans.

\begin{definition}[\Cref{def:tensor-reflexive}, \Cref{def:tg-marking}]
   We say that a complete marked fan $(\Sigma,u )$ inside a real vector space $V$ is \textit{$\otimes$-generating} if for every cone $\tau\in \Sigma \setminus \{0\}$, some affine hyperplane in $V$ contains the marking vectors $u_{\rho}$ when $\rho\in\tau(1)$, and has the marking vectors $u_{\rho}$, for the remaining rays of the cones containing $\tau$, strictly on one side. If the fan $\Sigma$ is in $N_{\bb R}$, then we call it \emph{$\otimes$-generating} if $(\Sigma,u)$ is $\otimes$-generating, where $u$ is given by primitive ray generators. 
\end{definition}

Our second main result is a derived category analogue of the Kleiman criterion for ampleness.

\begin{theorem}[\Cref{thm:curve-criterion}]\label{thm:B}
    Let $X$ be a proper toric variety over $k$ and $\ecal L$ a line bundle on $X$. Then $\ecal L$ is not $\otimes$-generating if and only if there is an irreducible complete curve $C \subseteq X$ with $\deg(\ecal L|_C) = 0$.
\end{theorem}

The curve $C$ is often \emph{not} torus-invariant. Indeed, on $\bb P^1 \times \bb P^1$, the line bundle $\ecal O(1,-1)$ has nonzero degree on the four invariant curves, but its restriction to the diagonal has degree $0$ (\Cref{ex:diagonal}). Producing such curves is the main work of \Cref{sec:curves}. The key step (\Cref{prop:degree-zero-curve}) is to construct a morphism from a projective line over a finite extension of $k$ using an integral relation among the ray generators, obtained from a variant of Farkas' lemma. The pullback of the divisor has degree zero and is therefore trivial. As a byproduct of the construction, we show that the closed cone of strongly movable curves of a complex projective toric variety, dual to the pseudoeffective cone by \cite{boucksom2013pseudo}*{Theorem 2.2}, is rational polyhedral. Every extremal ray is spanned by the class of an irreducible curve meeting the torus (\Cref{cor:sme-polyhedral}). The proof also implies that on a proper toric variety, a line bundle is $\otimes$-generating if and only if it is nowhere torsion in the sense of \cite{itoolander2026botamestack}, i.e., for any finite morphism $f: C \to X$ from a normal projective curve $C$, $f^*\ecal L$ is not torsion. 

Third, we study when $\otimes$-generating line bundles exist on complete toric varieties.

\begin{theorem}[\Cref{prop:existence-threefolds}, \Cref{thm:products}]\label{thm:C} \
    \begin{enumerate}
        \item Every simplicial proper toric variety of dimension $\leq 3$ admits a $\otimes$-generating line bundle.
        \item Let $X_1, \dots, X_n$ be positive-dimensional complete toric varieties with $n \geq 2$. The product $\prod_i X_i$ admits a $\otimes$-generating line bundle if and only if every $X_i$ is projective, in which case a line bundle is $\otimes$-generating if and only if it is ample or anti-ample. \qedhere 
    \end{enumerate}
\end{theorem}

For a non-projective complete toric variety $X$, the product $X \times \bb P^1$ therefore has no $\otimes$-generating line bundle. Taking $X$ smooth, non-projective, and three-dimensional produces smooth proper toric fourfolds with no $\otimes$-generating line bundle, and by part (i) the dimension $4$ is the smallest possible for smooth examples. Without smoothness, the phenomenon occurs already in dimension $3$ and Picard number $1$: a complete toric threefold of Fujino \cite{fujino2005kleiman}*{Lemma 2.2}, built to violate Kleiman's ampleness criterion, has every line bundle of degree $0$ on a fixed curve (\Cref{ex:fujino}).

As applications of these criteria, we construct several new examples of smooth proper varieties with $\otimes$-generating canonical bundles, to which the generalized Bondal--Orlov reconstruction applies. For example, we prove the following existence result, answering \cite{ito2025polarizations}*{Question 4.8}, affirmatively.  
\begin{theorem}[\Cref{thm:Xa,thm:Ya}]
    There exists an infinite family of smooth non-projective proper toric threefolds with $\otimes$-generating canonical bundle. 
\end{theorem}
More concretely, from the classification of smooth complete toric threefolds of Picard number $5$ by Oda and Miyake \cite{oda1978torus}*{Theorem 9.6}, the family $Z''_{13}(a,b,c,d)$, studied by Fujino and Sato \cite{fujino2024non}, contains non-projective members for $b \neq 0$, $a \neq c$. Within it, the threefolds $X_a := Z''_{13}(a,2,a+2,a+1)$ are smooth and non-projective with $\otimes$-generating canonical bundles for every $a \in \bb Z$, and $Y_a := Z''_{13}(a,2,a+2,a)$ has a $\otimes$-generating canonical bundle if and only if $a \neq -1, -2$ (\Cref{thm:Xa,thm:Ya}). 

We also construct examples with $\otimes$-generating canonical bundles, utilizing $\otimes$-generating fans and the derived analogue of Kleiman's criterion. Indeed, the iterated blow-up of $\bb P^3$ along a set $S$ of torus-invariant lines, taken in any order (\Cref{prop:P3-blowup-obstructions}), has $\otimes$-generating fan if and only if $S$ is the empty set, a single line, a pair of disjoint lines, a triple of concurrent lines, or all six lines. In addition, we prove that a projective bundle $\bb P_X(\ecal O_X^{\oplus n} \oplus \ecal L)$ over a proper Gorenstein toric variety $X$ has a $\otimes$-generating canonical bundle if $\omega_X \otimes \ecal L$ is anti-ample and $\omega_X \otimes \ecal L^{\otimes - n}$ is $\otimes$-generating. We study the case when $X$ is a Hirzebruch surface more closely and produce non-Fano examples that can still be reconstructed from their derived categories.

Finally, we study obstructions for $\otimes$-generation of canonical bundles. For convenience, we say an invariant subvariety $Z$ of a proper Gorenstein toric variety $X$ is \textit{obstructing} if $\omega_X|_Z$ is neither big nor anti-big. By \Cref{thm:A}, the existence of an obstructing invariant subvariety is indeed the obstruction for $\otimes$-generation of the canonical bundle, and hence we cannot apply the Bondal--Orlov reconstruction. In the section, we first characterize obstructing invariant curves in every dimension via wall relations (\Cref{prop:wall-criterion}) and classify such curves in dimension $3$ (\Cref{prop:wall-types-smooth}). Next, we observe that in certain cases, obstructing invariant divisors indeed contribute to non-standard autoequivalences of the derived category, via Horja's EZ spherical twist, obstructing conclusion (ii) of the reconstruction (\Cref{prop-obstructionproduceseztwist}).  


The rest of this paper is organized as follows. In \Cref{sec:conventions}, we fix notation, recall the results of \cite{ito2025derived} that we use, and characterize big line bundles on proper toric varieties together with their cone. In \Cref{sec:criterion}, we prove \Cref{thm:A} and treat first examples, including toric Fano varieties and Hirzebruch surfaces. We then introduce $\otimes$-generating fans. In \Cref{sec:curves}, we prove \Cref{thm:B} and the polyhedrality of the movable cone of curves. In \Cref{sec:existence}, we prove \Cref{thm:C} and present the non-existence examples.  In \Cref{sec:canonical}, we consider examples of $\otimes$-generating canonical bundles based on previous sections. In particular, we prove the results on the family $Z''_{13}$, the blow-ups of $\bb P^3$, and toric projective bundles. Many of these constitute new examples to which the Bondal--Orlov reconstruction applies. In \Cref{sec:obstruction}, we classify the obstructing walls for invariant curves and construct non-standard autoequivalences associated to certain obstructing divisors.

\paragraph{Acknowledgment} DI would like to thank John S. Nolan and Noah Olander for inspiring discussions. 

\paragraph{Code availability}\label{par:code-availability}
The SageMath \cite{sagemath} code for the bigness and $\otimes$-generation criteria and the computations in \Cref{thm:Xa,thm:Ya,prop:P3-blowup-obstructions} is available at 

\begin{center}
\texttt{\url{https://anonymous.4open.science/r/toric-tensor-generating-77D1/}}.
\end{center}

\paragraph{AI Usage Disclosure} 

The authors developed most mathematical contributions and initial drafts before March 2026, without AI assistance other than reference search. 

DI used ChatGPT 5.1 and 5.2 to (i) search for references, particularly, appropriate references in \cite{cox2011toric} for \Cref{sec:criterion}, and (ii) write Python code for \texttt{is\_tensor\_generating.py} and \texttt{Print\_PD\_for\_Xa\_Ya.py}, and other code for experiments and example searches. DI revised the first code by himself and is responsible for the correct implementation. DI also used ChatGPT 5.6 for reference searches. DI is responsible for the correctness of the outputs of the second code.

MRZ used ChatGPT 5.5 and 5.6 to search for references (in particular, in \Cref{subsec:Z13}) and to revise a few intermediate drafts. MRZ used ChatGPT 6 Sol to create \Cref{fig:P3-lines} and modify some existing plots. 

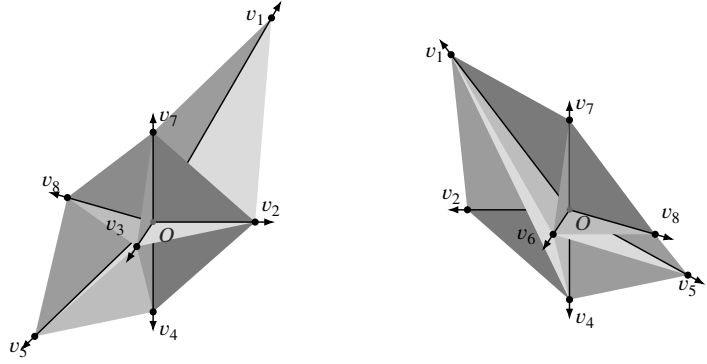
\begin{figure}[t]
    \centering



\begin{tikzpicture}[
    x={(1.25cm,0.00cm)},
    y={(0.20cm,0.30cm)},
    z={(0.00cm,1.10cm)},
    cls plot, scale=1.08]
    \let\wallfill\CLSwall
    \let\fanray\CLSray
    \colorlet{g1}{clsgray1}
    \colorlet{g2}{clsgray2}
    \colorlet{g3}{clsgray3}
    \colorlet{g4}{clsgray4}

    \begin{scope}
        \coordinate (O)  at ( 0, 0, 0);
        \coordinate (v1) at ( 1, 1, 2);
        \coordinate (v2) at ( 1, 0, 0);
        \coordinate (v3) at ( 0,-1, 0);
        \coordinate (v4) at ( 0, 0,-1);
        \coordinate (v5) at (-1,-1, -1);
        \coordinate (v6) at ( 0, 1, 0);
        \coordinate (v7) at ( 0, 0, 1);
        \coordinate (v8) at (-1, 1, 0);
    
    
        \wallfill{v5}{v8}{g3} 
        \wallfill{v7}{v8}{g3!50!g4} 
        \wallfill{v1}{v7}{g3}   
        \wallfill{v1}{v2}{g1}  
        \wallfill{v2}{v4}{g4}   
        \wallfill{v4}{v5}{g2}   
    
        \fanray{v1}
        
        \wallfill{v3}{v5}{g1} 
        \wallfill{v3}{v4}{g3}   
    
        \fanray{v4}
        \fanray{v5}
    
        \wallfill{v3}{v8}{g2} 
    
        \fanray{v8}
        
        \wallfill{v2}{v7}{g4}   
        \wallfill{v3}{v7}{g3} 
        \wallfill{v2}{v3}{g1}   
    
        \fanray{v2}
        \fanray{v3}
        \fanray{v7}
        
        \node[origin] at (O) {};
        \foreach \v in {v1,v2,v3,v4,v5,v7,v8} {
            \node[gen] at (\v) {};
        }
    
        \node[lbl, anchor=east] at (v1) {$v_1$};
        \node[lbl, anchor=south west] at (v2) {$v_2$};
        \node[lbl, anchor=south east, xshift=-2pt, yshift=2pt] at (v3) {$v_3$};
        \node[lbl, anchor=north west, yshift=-3pt] at (v4) {$v_4$};
        \node[lbl, anchor=north east] at (v5) {$v_5$};
        \node[lbl, anchor=south west] at (v7) {$v_7$};
        \node[lbl, anchor=south east] at (v8) {$v_8$};
        \node[lbl, text=gray!30!black, anchor=north west] at (O) {$O$};
    \end{scope}

    \begin{scope}[xshift=5.1cm, yshift=1.5mm]
        \coordinate (O)  at ( 0, 0, 0);
        \coordinate (v1) at (-1,-1, 2);
        \coordinate (v2) at (-1, 0, 0);
        \coordinate (v3) at ( 0, 1, 0);
        \coordinate (v4) at ( 0, 0,-1);
        \coordinate (v5) at ( 1, 1, -1);
        \coordinate (v6) at ( 0,-1, 0);
        \coordinate (v7) at ( 0, 0, 1);
        \coordinate (v8) at ( 1,-1, 0);
    

        \fanray{v3}
    
        \wallfill{v1}{v7}{g4}
        \wallfill{v8}{v5}{g3}
        \wallfill{v5}{v4}{g3}
        \wallfill{v4}{v2}{g4}
        \wallfill{v2}{v1}{g3}

        \fanray{v2}
        
        \wallfill{v1}{v4}{g1}

        \wallfill{v6}{v1}{g2}
        \wallfill{v6}{v4}{g2}

        \fanray{v1}
        \fanray{v4}

        \wallfill{v6}{v7}{g3}
        \wallfill{v6}{v5}{g1}

        \fanray{v5}
        \fanray{v7}

        \wallfill{v6}{v8}{g2}
        \wallfill{v7}{v8}{g4}

        \fanray{v6}
        \fanray{v8}
        
        \node[origin] at (O) {};
        \foreach \v in {v1,v6,v2,v4,v5,v7,v8} {
            \node[gen] at (\v) {};
        }
    
        \node[lbl, anchor=east] at (v1) {$v_1$};
        \node[lbl, anchor=south east] at (v2) {$v_2$};
        \node[lbl, anchor=north west, yshift=-3pt] at (v4) {$v_4$};
        \node[lbl, anchor=north] at (v5) {$v_5$};
        \node[lbl, anchor=east, xshift=-4pt, yshift=-1pt] at (v6) {$v_6$};
        \node[lbl, anchor=south west] at (v7) {$v_7$};
        \node[lbl, anchor=south west, yshift=2pt] at (v8) {$v_8$};
        \node[lbl, text=gray!30!black, anchor=north west] at (O) {$O$};
    
    \end{scope}
\end{tikzpicture}


    \caption{The fan of $Y_{-1}=Z_{13}''(-1,2,1,-1)$ as in \Cref{thm:Ya}. Since $v_2,v_3,v_4$ and $v_5$ are coplanar, $\omega_{Y_{-1}}$ is not $\otimes$-generating.} \label{fig:Y-1}
\end{figure}

XZ used ChatGPT 5.5 and 5.6 to search for references (in particular, in \Cref{subsec:products}) and to help prepare \Cref{fig:X-1}, \Cref{fig:Y-1}, and \Cref{fig:bundle-over-F4}.

Additionally, this paper went through checks by ChatGPT 6 Astra during the final drafting phase and we corrected some gaps and typographical errors, as well as added clarifications based on its suggestions, without changing main results.  

We reviewed and verified the AI-assisted material incorporated into the paper and take full responsibility for the contents of the paper.

\section{Preliminaries}\label{sec:conventions}

Throughout the paper, $k$ denotes a field. A variety is an integral separated scheme of finite type over $k$. Triangulated categories and exact functors are assumed to be $k$-linear. For toric varieties, we follow the conventions of \cite{cox2011toric}. The results we cite from \cite{cox2011toric} are stated there over $\bb C$, although their proofs are combinatorial and hold over $k$ verbatim (see \cite{danilov1978geometry}). For polytopes, we follow \cite{beck2007computing}*{\S 2.1}.

\subsection{Toric conventions}\label{subsec:toric-conventions}
Let $N \iso \bb Z^n$ be a lattice and $M := \operatorname{Hom}_{\bb Z}(N, \bb Z)$ its dual, with pairing $\bra{-,-}\colon M \times N \to \bb Z$. For a fan $\Sigma$ in $N_{\bb R} := N \otimes_{\bb Z} \bb R$, we write $X = X_\Sigma$ for the associated toric variety with dense torus $T_N$. The variety $X_\Sigma$ is proper if and only if $\Sigma$ is complete, i.e., $|\Sigma| = N_{\bb R}$ \cite{cox2011toric}*{Theorem 3.4.6}.

We write $\Sigma(r)$ for the set of $r$-dimensional cones of $\Sigma$ and $\sigma(1)$ for the set of rays of a cone $\sigma \in \Sigma$. The relation $\tau \prec \sigma$ means that $\tau$ is a face of $\sigma$. Each ray $\rho \in \Sigma(1)$ has a primitive generator $u_\rho \in N$ and determines a $T_N$-invariant prime Weil divisor $D_\rho \subset X$.

For a cone $\tau \in \Sigma$, we write $V(\tau) \subseteq X$ for the closure of the corresponding torus orbit \cite{cox2011toric}*{\S 3.2}, so that $\dim V(\tau) = n - \dim \tau$ and $V(\rho) = D_\rho$ for a ray $\rho$. Set $N(\tau) := N/(\bra{\tau} \cap N)$, where $\bra{\tau}$ denotes the linear span, and let $p_\tau \colon N \to N(\tau)$ be the projection. By \cite{cox2011toric}*{Proposition 3.2.7}, $V(\tau)$ is the toric variety associated with the \emph{star fan}
\[
\operatorname{Star}(\tau) := \{\bar{\sigma} := p_\tau(\sigma) \subset N(\tau)_{\bb R} \mid \tau \prec \sigma \in \Sigma\}
\]
in $N(\tau)_{\bb R}$. We use the following observation repeatedly: for cones $\tau \prec \sigma$ in $\Sigma$, the orbit closure $V_{\operatorname{Star}(\tau)}(\bar \sigma) \subseteq V(\tau)$ of the cone $\bar\sigma \in \operatorname{Star}(\tau)$ coincides with the orbit closure $V(\sigma) \subseteq X$. In particular, every $T_{N(\tau)}$-invariant prime divisor of $V(\tau)$ is an orbit closure of $X$.

A $T_N$-invariant Weil divisor on $X$ is an integral linear combination $D = \sum_{\rho \in \Sigma(1)} a_\rho D_\rho$. The group $\operatorname{Div}_{T_N}(X)$ of $T_N$-invariant Cartier divisors fits into an exact sequence
\[
M \to \operatorname{Div}_{T_N}(X) \to \operatorname{Pic}(X) \to 0,
\]
where $m \in M$ maps to $\operatorname{div}(\chi^m) = \sum_\rho \bra{m, u_\rho} D_\rho$ \cite{cox2011toric}*{Theorem 4.2.1}. In particular, every line bundle on $X$ is isomorphic to $\ecal O_X(D)$ for some $D \in \operatorname{Div}_{T_N}(X)$. A $T_N$-invariant Cartier divisor $D$ is specified by its \emph{Cartier data} $\{m_\sigma \in M\}_{\sigma \in \Sigma_{\max}}$, indexed by the maximal cones and characterized by $\bra{m_\sigma, u_\rho} = -a_\rho$ for all $\rho \in \sigma(1)$ \cite{cox2011toric}*{Theorem 4.2.8}. We extend the data to all cones by setting $m_\tau := m_\sigma$ for a chosen maximal cone $\sigma \supseteq \tau$. Then $\bra{m_\tau, u_\rho} = -a_\rho$ for all $\rho \in \tau(1)$, which is the only property we use. A divisor $D$ is $\bb Q$-Cartier if $lD$ is Cartier for some integer $l > 0$.

\subsection{Derived categories and \texorpdfstring{$\otimes$-generation}{tensor-generation}}\label{subsec:derived-conventions}
For a scheme $X$, we write $\perf(X)$ for the triangulated category of perfect complexes. A set $S$ of objects of $\perf(X)$ \emph{classically generates} $\perf(X)$ if the smallest thick subcategory of $\perf(X)$ containing $S$ is $\perf(X)$ itself \cite{bondal2002generators}. Following \cite{ito2025derived}*{Definition 2.1}, where the notion of \cite{ito2025polarizations}, introduced there under the name $\tens$-ample, receives its present name, we say a line bundle $\ecal L$ on $X$ is \emph{$\otimes$-generating} if the set of powers $\{\ecal L^{\tens m} \mid m \in \bb Z\}$ classically generates $\perf(X)$. Since this set coincides with the set of powers of $\ecal L^{-1}$, a line bundle is $\otimes$-generating if and only if its inverse is.

A line bundle $\ecal L$ on a proper variety $X$ of dimension $n$ is \emph{big} if $\limsup_{l \to \infty} \dim_k \h^0(X, \ecal L^{\tens l})/l^n > 0$, and \emph{anti-big} if $\ecal L^{-1}$ is big. On proper varieties, this notion of bigness agrees with the one used in \cite{ito2025derived} by \cite{ito2025derived}*{Lemma 2.9}. For a $\bb Q$-Cartier divisor $D$, we say $D$ is big (resp.\ anti-big) if $\ecal O_X(lD)$ is big (resp.\ anti-big) for one, equivalently any, integer $l > 0$ such that $lD$ is Cartier.

Three results of Ito--Olander enter our arguments. First, their derived Nakai--Moishezon criterion \cite{ito2025derived}*{Theorem 1.1}: a line bundle $\ecal L$ on a proper $k$-scheme $X$ is $\otimes$-generating if and only if the restriction $\ecal L|_Z$ is big or anti-big for every closed subvariety $Z \subseteq X$. Here a closed subvariety is an integral closed subscheme, and $X$ itself counts among them when $X$ is a variety. Second, a criterion of Ito--Olander \cite{ito2025derived}*{Lemma 3.19}, which we state in the special case of a single section. Let $\ecal L$ be a line bundle on a noetherian scheme $X$ and $s \in \h^0(X, \ecal L^{\tens l})$ a section such that $X_s$ is quasi-affine. If the restriction of $\ecal L$ to the reduced zero scheme $Z(s)$ of $s$ is $\otimes$-generating, then $\ecal L$ is $\otimes$-generating. Together with the invariance of $\otimes$-generation under passing to irreducible components \cite{ito2025derived}*{Corollary 1.3}, this criterion reduces $\otimes$-generation to the components of $Z(s)$. Third, the case of curves \cite{ito2025derived}*{Proposition 5.1}: a line bundle on an integral projective curve over $k$ is $\otimes$-generating if and only if its degree is nonzero.

\subsection{Big line bundles on proper toric varieties}\label{subsec:bigness}

We now review how to characterize big line bundles on a proper toric variety by the dimensions of the associated polytopes and describe the big cone. Let $X = X_\Sigma$ be a proper toric variety of dimension $n$ and let $D = \sum_\rho a_\rho D_\rho$ be a $T_N$-invariant Cartier divisor. Recall from \cite{cox2011toric}*{Proposition 4.3.3} that
\begin{equation}\label{eq:sections-lattice-points}
\h^0(X, \ecal O_X(D)) = \bigoplus_{m \in P_D \cap M} k \cdot \chi^m,
\qquad
P_D := \{m \in M_{\bb R} \mid \bra{m, u_\rho} \geq -a_\rho \ \text{for all } \rho \in \Sigma(1)\}.
\end{equation}
Since $\Sigma$ is complete, $P_D$ is a rational convex polytope \cite{cox2011toric}*{Proposition 4.3.8}, and $P_{lD} = lP_D$ for every integer $l \geq 0$ \cite{cox2011toric}*{p.~192}. The \emph{support function} $\varphi_D \colon N_{\bb R} \to \bb R$ is the function that is linear on each cone of $\Sigma$ and satisfies $\varphi_D(u_\rho) = -a_\rho$ for every ray $\rho$. The support function determines and is determined by $D$ \cite{cox2011toric}*{Theorem 4.2.12}. The following seems to be well-known to experts.

\begin{lemma}\label{lem:toric-bigness}
    Let $X = X_\Sigma$ be a proper toric variety of dimension $n$ and let $\ecal L$ be a line bundle on $X$. The following are equivalent.
    \begin{enumerate}
        \item The line bundle $\ecal L$ is big.
        \item For one, equivalently every, $T_N$-invariant Cartier divisor $D$ with $\ecal L \iso \ecal O_X(D)$, the polytope $P_D$ satisfies $\dim P_D = n$.
        \item For one, equivalently every, $T_N$-invariant Cartier divisor $D$ with $\ecal L \iso \ecal O_X(D)$, there is a linear function $f \colon N_{\bb R} \to \bb R$ with $f > \varphi_D$ on $N_{\bb R} \setminus \{0\}$.
        \item There exist an integer $l > 0$ and a section $s \in \h^0(X, \ecal L^{\tens l})$ with $X_s = T_N$. \qedhere
    \end{enumerate}
\end{lemma}
\begin{proof}
    By the exact sequence in \Cref{subsec:toric-conventions}, $\ecal L \iso \ecal O_X(D)$ for some $T_N$-invariant Cartier divisor $D$, and any two such divisors differ by $\operatorname{div}(\chi^m)$ for some $m \in M$, which translates $P_D$ by $-m$ and changes $\varphi_D$ by the linear function $-m$. Conditions (ii) and (iii) are therefore independent of the choice of $D$, and it suffices to prove the equivalences for a fixed $D$.

    (ii) $\iff$ (iii). A linear function $f$ satisfies $f > \varphi_D$ on $N_{\bb R} \setminus \{0\}$ if and only if $f(u_\rho) > -a_\rho$ for every ray $\rho$, since both $f$ and $\varphi_D$ are linear on each cone of the complete fan $\Sigma$ and the rays of a cone span it. The latter condition says that $f$ is an interior point of $P_D$, and a convex polytope in $M_{\bb R}$ has nonempty interior if and only if it has dimension $n$.

    (i) $\iff$ (ii). By \Cref{eq:sections-lattice-points}, $\dim_k \h^0(X, \ecal O_X(lD)) = \#(lP_D \cap M)$ for $l > 0$. Choose an integer $m_0 > 0$ such that $m_0 P_D$ is a lattice polytope. By Ehrhart's theorem \cite{beck2007computing}*{Theorem 3.8}, $\#(l m_0 P_D \cap M)$ is a polynomial in $l$ of degree $\dim P_D$. If $\dim P_D = n$, then $\dim_k \h^0(X, \ecal O_X(l m_0 D))$ grows like $l^n$, so $\ecal L$ is big. If $\dim P_D < n$, then $P_D$ is contained in a polytope $Q$ of dimension $d := \dim P_D$ that is a translate of a lattice polytope by a rational vector, and $\#(lP_D \cap M) \leq \#(lQ \cap M) = O(l^{d})$, so $\ecal L$ is not big.

    (ii) $\implies$ (iv). Since $\dim P_D = n$, the interior of $lP_D = P_{lD}$ contains a lattice point $m$ for $l$ sufficiently large. Let $s_m \in \h^0(X, \ecal O_X(lD))$ be the corresponding section. Then
    \[
    \operatorname{div}(s_m) = \operatorname{div}(\chi^m) + lD = \sum_\rho (\bra{m, u_\rho} + l a_\rho) D_\rho,
    \]
    and $\bra{m, u_\rho} + l a_\rho > 0$ for every ray $\rho$ because $m$ is interior to $P_{lD}$. Hence $X_{s_m} = X \setminus \bigcup_{\rho} D_\rho = T_N$.

    (iv) $\implies$ (i). This follows by \cite{ito2025derived}*{Lemmas 2.8 and 2.9} since $T_N$ is affine. 
\end{proof}

\begin{remark}\label{rem:Q-Cartier-bigness}
    \Cref{lem:toric-bigness} extends to $\bb Q$-Cartier divisors. For a $T_N$-invariant $\bb Q$-Cartier divisor $D$, choose $l > 0$ with $lD$ Cartier. Then $D$ is big if and only if $\dim P_D = n$, since $P_{lD} = l P_D$. Similarly, condition (iii) holds for $D$ if and only if it holds for $lD$.
\end{remark}

\begin{remark}\label{rem:known-cases}
    When $\ecal L$ is nef, the equivalence of (i) and (ii) in \Cref{lem:toric-bigness} is shown in \cite{cox2011toric}*{Lemma 9.3.9}. When $X$ is simplicial and semiprojective, the characterization also follows from \cite{cox2011toric}*{Lemma 15.1.8}, which states that $\operatorname{PsEff}(X) = \sum_{\rho \in \Sigma(1)} \bb R_{\geq 0}\,[D_\rho]$ in $N^1(X)_{\bb R}$. Our result gives the following generalization.
\end{remark}

\begin{corollary}\label{cor:big-cone}
    Let $X = X_\Sigma$ be a proper toric variety.
    \begin{enumerate}
        \item A $T_N$-invariant $\bb Q$-Cartier divisor $D$ is big if and only if $D \sim_{\bb Q} \sum_\rho a_\rho D_\rho$ with $a_\rho \in \bb Q_{>0}$ for every ray $\rho$.
        \item The variety $X$ carries a big line bundle if and only if there are integers $a_\rho > 0$ for $\rho \in \Sigma(1)$ such that $\sum_\rho a_\rho D_\rho$ is Cartier. \qedhere 
    \end{enumerate}
\end{corollary}
\begin{proof}
    (i) If $D \sim_{\bb Q} \sum_\rho a_\rho D_\rho$ with all $a_\rho > 0$, then $P_{\sum a_\rho D_\rho}$ contains a ball around the origin, so $D$ is big by \Cref{lem:toric-bigness} and \Cref{rem:Q-Cartier-bigness}. Conversely, if $D$ is big, take $l > 0$ with $lD$ Cartier. By \Cref{lem:toric-bigness}(iv) applied to $\ecal O_X(lD)$, there are $l' > 0$ and $s \in \h^0(X, \ecal O_X(l'lD))$ with $X_s = T_N$. Then $\operatorname{div}(s) = \sum_\rho b_\rho D_\rho$ with every $b_\rho > 0$, and $D \sim_{\bb Q} \frac{1}{l'l} \operatorname{div}(s)$.

    (ii) If such $a_\rho$ exist, $\ecal O_X(\sum_\rho a_\rho D_\rho)$ is big as in (i). Conversely, a big line bundle is isomorphic to $\ecal O_X(D)$ for an invariant Cartier divisor $D$, and by (i) $D\sim_{\bb Q} \sum a_\rho D_\rho$, $a_\rho \in \bb Q_{>0}$. This divisor is $\bb Q$-Cartier, and clearing denominators by a positive integer makes it Cartier with positive integral coefficients.
\end{proof}

\begin{remark}\label{rem:big-cone-simplicial}
    Suppose $X$ is in addition simplicial, so that every $T_N$-invariant $\bb R$-divisor is $\bb R$-Cartier and $N^1(X)_{\bb R}$ is the quotient of $\bigoplus_\rho \bb R \, [D_\rho]$ by the relations $\sum_\rho \bra{m, u_\rho}[D_\rho] = 0$ for $m \in M$ \cite{cox2011toric}*{Theorem 4.1.3, Proposition 4.2.7, and Proposition 6.3.15}. \Cref{cor:big-cone} then identifies the big cone:
    \[
    \operatorname{Big}(X) = \Bigl\{\textstyle\sum_\rho a_\rho [D_\rho] \in N^1(X)_{\bb R} \Bigm| a_\rho > 0 \ \text{for all } \rho \in \Sigma(1)\Bigr\}. \qedhere 
    \]
\end{remark}

\section{The combinatorial criterion}\label{sec:criterion}

In this section, we prove \Cref{thm:A}. We first reduce $\otimes$-generation to bigness of restrictions to orbit closures, then compute these restrictions by polytopes, and finally reformulate the criterion as a hyperplane separation condition on markings of the fan.

\subsection{The criterion on invariant subvarieties}\label{subsec:criterion-invariant-subvar}

\begin{theorem}\label{thm:orbit-closure}
    Let $X = X_\Sigma$ be a proper toric variety and $\ecal L$ a big line bundle on $X$. The following are equivalent.
    \begin{enumerate}
        \item The line bundle $\ecal L$ is $\otimes$-generating.
        \item The restriction $\ecal L|_{D_\rho}$ is $\otimes$-generating for every ray $\rho \in \Sigma(1)$.
        \item The restriction $\ecal L|_{V(\tau)}$ is big or anti-big for every cone $\tau \in \Sigma$. \qedhere
    \end{enumerate}
\end{theorem}
\begin{proof}
    (ii) $\implies$ (i). By \Cref{lem:toric-bigness}, there are $l > 0$ and $s \in \h^0(X, \ecal L^{\tens l})$ with $X_s = T_N$, so that the reduced zero scheme of $s$ satisfies $Z(s) = \bigcup_{\rho \in \Sigma(1)} D_\rho$ with irreducible components $D_\rho$. If every $\ecal L|_{D_\rho}$ is $\otimes$-generating, then $\ecal L|_{Z(s)}$ is $\otimes$-generating by \cite{ito2025derived}*{Corollary 1.3}, and since $X_s = T_N$ is affine, $\ecal L$ is $\otimes$-generating by \cite{ito2025derived}*{Lemma 3.19}.

    (i) $\implies$ (iii). Each orbit closure $V(\tau)$ is a closed subvariety of $X$, so this is part of \cite{ito2025derived}*{Theorem 1.1}.

 (iii) $\implies$ (ii). We argue by induction on $\dim X$ that any line bundle satisfying (iii) is $\otimes$-generating. The case $\dim X=0$ is immediate. Condition (iii) at $\tau=\{0\}$ makes $\ecal L$ big or anti-big, so after replacing $\ecal L$ by $\ecal L^{-1}$ if necessary, we may assume it is big. For each ray $\rho$, every orbit closure in $D_\rho$ is an orbit closure in $X$ by the observation in \Cref{subsec:toric-conventions}. Thus $\ecal L|_{D_\rho}$ satisfies (iii) and is $\otimes$-generating by the inductive hypothesis. 
\end{proof}

\Cref{thm:orbit-closure} reduces $\otimes$-generation to bigness of the restrictions $\ecal L|_{V(\tau)}$. The following allows a completely combinatorial description of $\otimes$-generation. For a $T_N$-invariant Cartier divisor $D = \sum_\rho a_\rho D_\rho$ and a cone $\tau \in \Sigma$, define the \emph{restriction polytope}
\begin{equation}\label{eq:restriction-polytope}
P_D(\tau) := \{m \in M_{\bb R} \mid \bra{m, u_\rho} = -a_\rho \ \forall \rho \in \tau(1), \ \text{and} \ \bra{m, u_\rho} \geq -a_\rho \ \forall \rho \in \sigma(1), \ \forall \sigma \succ \tau\}.
\end{equation}
For $\tau = \{0\}$, the equalities are vacuous and every cone contains $\tau$, so $P_D(\{0\}) = P_D$.

\begin{lemma}\label{lem:restriction-polytope}
    Let $X = X_\Sigma$ be a proper toric variety and $D$ a $T_N$-invariant Cartier divisor on $X$, with restriction polytopes as in \Cref{eq:restriction-polytope}. For every cone $\tau \in \Sigma$, there is a $T_{N(\tau)}$-invariant Cartier divisor $D_\tau$ on $V(\tau)$ such that
    \[
    \ecal O_X(D)|_{V(\tau)} \iso \ecal O_{V(\tau)}(D_\tau)
    \qquad \text{and} \qquad
    P_{D_\tau} = P_D(\tau) - m_\tau
    \]
    for any Cartier data $\{m_\sigma\}_{\sigma \in \Sigma}$ of $D$. In particular, $\dim P_{D_\tau} = \dim P_D(\tau)$.
\end{lemma}
\begin{proof}
    Fix Cartier data $\{m_\sigma\}_{\sigma \in \Sigma}$ of $D$, so that $\bra{m_\sigma, u_\rho} = -a_\rho$ for all $\rho \in \sigma(1)$. The divisor $\tilde{D}_\tau := D + \operatorname{div}(\chi^{m_\tau})$ has Cartier data $\{m_\sigma - m_\tau\}_{\sigma \in \Sigma}$ \cite{cox2011toric}*{Exercise 4.2.6} and satisfies $\ecal O_X(D) \iso \ecal O_X(\tilde D_\tau)$. For $\sigma \succ \tau$, the linear functions $m_\sigma$ and $m_\tau$ agree on $\tau$, since both compute the support function $\varphi_D$ there. Hence
    \[
    m_{\bar\sigma} := m_\sigma - m_\tau \in M \cap \tau^\perp = N(\tau)^\vee \qquad (\sigma \succ \tau),
    \]
    and the collection $\{m_{\bar\sigma}\}_{\bar\sigma \in \operatorname{Star}(\tau)}$ is Cartier data of a $T_{N(\tau)}$-invariant Cartier divisor $D_\tau$ on $V(\tau)$. The restrictions of $\tilde D_\tau$ and $D_\tau$ have the same transition functions on the invariant affine cover of $V(\tau)$, so
    \[
    \ecal O_X(D)|_{V(\tau)} \iso \ecal O_X(\tilde D_\tau)|_{V(\tau)} \iso \ecal O_{V(\tau)}(D_\tau),
    \]
    as in the proof of \cite{cox2011toric}*{Proposition 6.3.8}.

    We now compute $P_{D_\tau}$. Write $\bra{-,-}'$ for the pairing between $M \cap \tau^\perp$ and $N(\tau)$, and $p := p_\tau \colon N \to N(\tau)$ for the projection. By definition of Cartier data, a point $m \in (M \cap \tau^\perp)_{\bb R}$ lies in $P_{D_\tau}$ if and only if
    \begin{equation}\label{eq:PDtau-star}
    \bra{m - m_{\bar\sigma}, v}' \geq 0
    \qquad \text{for every } \bar\sigma \in \operatorname{Star}(\tau) \text{ and every } v \in \bar\sigma.
    \end{equation}
    To see this, note that the divisor coefficient of $D_\tau$ along a ray $\bar\rho$ of $\operatorname{Star}(\tau)$ is $-\bra{m_{\bar\sigma}, u_{\bar\rho}}'$ for any $\bar\sigma$ containing $\bar\rho$. The defining inequalities of $P_{D_\tau}$ are therefore \Cref{eq:PDtau-star} for $v$ ranging over ray generators, and the general case follows because each $\bar\sigma$ is generated by its rays. On the other hand, for $\sigma \succ \tau$ and $\rho \in \sigma(1)$, we have $m - m_{\bar\sigma} \in \tau^\perp$ and $p(u_\rho) \in \bar\sigma$, whence
    \[
    \bra{m + m_\tau, u_\rho} \geq -a_\rho
    \iff \bra{m + m_\tau - m_\sigma, u_\rho} \geq 0
    \iff \bra{m - m_{\bar\sigma}, p(u_\rho)}' \geq 0.
    \]
    Since every cone $\bar\sigma \in \operatorname{Star}(\tau)$ is generated by the vectors $p(u_\rho)$ with $\rho \in \sigma(1)$, condition \Cref{eq:PDtau-star} holds for $m$ if and only if $\bra{m + m_\tau, u_\rho} \geq -a_\rho$ for all $\rho \in \sigma(1)$ and all $\sigma \succ \tau$. The equalities $\bra{m + m_\tau, u_\rho} = -a_\rho$ for $\rho \in \tau(1)$ hold automatically for $m \in \tau^\perp$, and conversely any element of $P_D(\tau)$ lies in $m_\tau + \tau^\perp$. Hence $P_{D_\tau} = P_D(\tau) - m_\tau$.
\end{proof}

Combining \Cref{thm:orbit-closure} with \Cref{lem:toric-bigness,lem:restriction-polytope} yields the polytope form of the criterion in \Cref{thm:A}.

\begin{theorem}\label{thm:implementable}
    Let $X = X_\Sigma$ be a proper toric variety and $D$ a $T_N$-invariant $\bb Q$-Cartier divisor on $X$. Then $\ecal O_X(lD)$ is $\otimes$-generating for $l > 0$ with $lD$ Cartier if and only if for every cone $\tau \in \Sigma$, either $\dim P_D(\tau)$ or $\dim P_{-D}(\tau)$ equals $\dim X - \dim \tau$.
\end{theorem}
\begin{proof}
    Since $P_{lD}(\tau) = l \, P_D(\tau)$, we may assume $D$ is Cartier and $l = 1$. Applying the case $\tau = \{0\}$ of the criterion, together with \Cref{lem:toric-bigness} and the identity $P_D(\{0\}) = P_D$, we may assume $\ecal O_X(D)$ is big or anti-big, and after replacing $D$ by $-D$, big. By \Cref{thm:orbit-closure}, $\ecal O_X(D)$ is $\otimes$-generating if and only if $\ecal O_X(D)|_{V(\tau)}$ is big or anti-big for every $\tau \in \Sigma$. By \Cref{lem:restriction-polytope}, $\ecal O_X(\pm D)|_{V(\tau)} \iso \ecal O_{V(\tau)}((\pm D)_\tau)$, and by \Cref{lem:toric-bigness} applied on the $(\dim X - \dim\tau)$-dimensional proper toric variety $V(\tau)$, the restriction is big if and only if $\dim P_{(\pm D)_\tau} = \dim P_{\pm D}(\tau)$ equals $\dim X - \dim \tau$.
\end{proof}

\begin{remark}\label{rem:intermediate-cones}
    In \Cref{thm:implementable}, the conditions at $\tau = \{0\}$ and at maximal cones need no verification once $D$ is big or anti-big. At $\tau = \{0\}$, the condition is bigness or anti-bigness of $D$ itself. At a maximal cone $\sigma$, the orbit closure $V(\sigma)$ is a point, and $P_D(\sigma) = \{m \in M_{\bb R} \mid \bra{m, u_\rho} = -a_\rho \ \forall \rho \in \sigma(1)\}$ contains $m_\sigma$ and is cut out by equations whose normal vectors $u_\rho$, $\rho \in \sigma(1)$, span $N_{\bb R}$, so $\dim P_D(\sigma) = 0$ always holds. For a big or anti-big $D$, it therefore suffices to consider the cones $\tau$ with $0 < \dim \tau < \dim X$.
\end{remark}

\subsection{Markings and the hyperplane criterion} \label{subsec:hyperplane-criterion}

We now give a description, in terms of the geometry of the fan, of when a big line bundle is $\otimes$-generating. Following the terminology of Ross and collaborators (\cite{nathanson2022tropical}*{\S 2.1}, \cite{nowak2023mixed}*{\S 2.1}), a \emph{marking} of a fan $\Sigma$ is a choice $u = (u_\rho)_{\rho \in \Sigma(1)}$ of a nonzero point $u_\rho$ on each ray $\rho$. For a $T_N$-invariant divisor $D = \sum_\rho a_\rho D_\rho$ with all $a_\rho > 0$, the \emph{marking associated with $D$} is
\[
u_D := (u_{D, \rho})_{\rho \in \Sigma(1)}, \qquad u_{D, \rho} := \frac{u_\rho}{a_\rho}.
\]

\begin{prop}\label{prop:hyperplane-test}
    Let $X = X_\Sigma$ be a proper toric variety and $\ecal L$ a big line bundle on $X$. Then there exist an integer $l > 0$ and a $T_N$-invariant $\bb Q$-Cartier divisor $D = \sum_\rho a_\rho D_\rho$ with all $a_\rho > 0$, $lD$ Cartier, and $\ecal L^{\tens l} \iso \ecal O_X(lD)$. For any such $D$, the line bundle $\ecal L$ is $\otimes$-generating if and only if for every cone $\tau \in \Sigma \setminus \{0\}$, there is an affine hyperplane $H \subset N_{\bb R}$ such that
    \begin{enumerate}
        \item $u_{D, \rho} \in H$ for every $\rho \in \tau(1)$, and
        \item the points $u_{D,\rho}$ for $\rho \in \sigma(1) \setminus \tau(1)$ with $\sigma \succ \tau$ all lie strictly on one side of $H$. \qedhere
    \end{enumerate}
\end{prop}
\begin{proof}
    According to \Cref{lem:toric-bigness}, for any invariant Cartier divisor $D$ representing $\ecal L$ there exists $m\in M_{\bb R}$ such that $a_\rho+\bra{m,u_\rho}>0$. We can perturb $m$ so that $m=m'/l$ with $l\in \bb N_{>0}, m'\in M$, and we can replace $D$ by $D+\mathrm{div}(\chi^{m'})/l$. Now let $D$ and $l$ be any choice as in the statement. For nonzero maximal cones, the required hyperplane is given by the Cartier data of $lD$. Since $\otimes$-generation is unchanged by taking a positive power and scaling would not change the dimension of $P_D(\tau)$, by \Cref{rem:intermediate-cones}, it suffices to show that for any cone $\tau$ with $0<\dim\tau<\dim X$, the condition above is equivalent to the condition in \Cref{thm:implementable}.

    If $P_D(\tau)$ has the expected dimension, then it must span the affine space
    \[
        A=\{m\in M_{\bb R}\mid \bra{m,u_{D,\rho}}=\frac{1}{a_\rho}\bra{m,u_\rho}=-1\ \ \forall \rho\in\tau(1)\}.
    \]
    The affine space $A$ is a translate of $\tau^\perp$, since $D$ is $\bb Q$-Cartier. Take $m$ to be in the interior of $P_{D}(\tau)$ inside $A$. Then for $\rho\in \sigma(1)\setminus\tau(1)$\ with\ $\sigma\succ\tau$, we must have $\bra{m,u_{D,\rho}}>-1$. Indeed, the corresponding functionals are nonconstant on $A$, since $u_\rho\notin\operatorname{span}(\tau)$ for such $\rho$, as $\sigma\cap\operatorname{span}(\tau)=\tau$. Then $H=\{u\in N_{\bb R}\mid\bra{m,u}=-1\}$ would satisfy the condition above. If $P_{-D}(\tau)$ has the expected dimension, we can use the same argument by taking $H=\{u\in N_{\bb R}\mid \bra{m, u}=1\}$.

    On the other hand, if there exists an $H$ satisfying the above condition, then the origin is not contained in $H$, or otherwise $H$ would induce a linear hyperplane in $N(\tau)_{\bb R}$ with all ray generators of $\operatorname{Star}(\tau)$ strictly on one side, which cannot happen since $\operatorname{Star}(\tau)$ is complete. If the origin lies in the same half-space as those $u_{D,\rho}$, then the corresponding closed half-space can be described as $\{u\in N_{\bb R}\mid \bra{m, u}\geq -1\}$ for some $m\in M_{\bb R}$. This would imply that $P_D(\tau)$ has $m$ inside its interior when viewed as a subset of $A$, hence having the expected dimension. If the origin lies on the other side, then the corresponding closed half-space can be described as $\{u\in N_{\bb R}\mid \bra{m, u}\geq 1\}$, and the same argument in $-A$ shows that $P_{-D}(\tau)$ has the expected dimension.
\end{proof}

\begin{corollary}[The hyperplane criterion]\label{cor:conv-aff}
    In the setting of \Cref{prop:hyperplane-test}, $\ecal L$ is $\otimes$-generating if and only if for every cone $\tau \in \Sigma \setminus \{0\}$,
    \[
    \operatorname{conv}\bigl(u_{D, \rho} \bigm| \rho \in \sigma(1) \setminus \tau(1), \ \sigma \succ \tau \bigr)
    \ \cap \
    \operatorname{aff}\bigl(u_{D, \rho} \bigm| \rho \in \tau(1)\bigr) = \emp.
    \]
    In particular, if $X$ is $\bb Q$-Gorenstein and $l > 0$ is an integer with $lK_X$ Cartier, then $\ecal O_X(lK_X)$ is $\otimes$-generating if and only if the above condition holds for the marking by primitive generators $u_\rho$.
\end{corollary}
\begin{proof}
    Write $S$ for the finite set of points $u_{D,\rho}$ with $\rho \in \sigma(1)\setminus\tau(1)$, $\sigma \succ \tau$, and $T$ for the set of points $u_{D,\rho}$ with $\rho \in \tau(1)$. It suffices to show that a hyperplane as in \Cref{prop:hyperplane-test} exists if and only if $\operatorname{conv}(S) \cap \operatorname{aff}(T) = \emp$. When $\tau$ is maximal, $S$ is empty, and the claim holds trivially. Otherwise, if $H$ is such a hyperplane, then $\operatorname{aff}(T) \subseteq H$ while $\operatorname{conv}(S)$ lies strictly on one side of $H$, so the two sets are disjoint. Conversely, suppose they are disjoint. The set $\operatorname{conv}(S)$ is compact convex and $\operatorname{aff}(T)$ is closed and convex, so there is a linear functional $m \in M_{\bb R}$ with
    \[
    \sup_{u \in \operatorname{aff}(T)} \bra{m, u} < \min_{u \in \operatorname{conv}(S)} \bra{m, u}.
    \]
    A linear functional bounded above on an affine subspace is constant on it, say with value $c$, and $H := \{u \mid \bra{m, u} = c\}$ contains $\operatorname{aff}(T) \supseteq T$ and has $S$ strictly on one side.

    For the second claim, we have $K_X = -\sum_\rho D_\rho$ \cite{cox2011toric}*{Theorem 8.2.3}, so $D := -lK_X = \sum_\rho l D_\rho$ is Cartier with all coefficients positive, hence big by \Cref{cor:big-cone}, and $\ecal L := \ecal O_X(-lK_X)$ satisfies the hypotheses of \Cref{prop:hyperplane-test} with this $D$ and marking $u_{D, \rho} = u_\rho / l$. Rescaling by $l$ carries the affine hull and the convex hull of the marking vectors to those of the primitive generators and preserves disjointness, so the desired condition for the marking holds if and only if it holds for the $u_\rho$. Since a line bundle is $\otimes$-generating if and only if its inverse is, the claim follows.
\end{proof}

\subsection{First examples}\label{subsec:first-examples}

\begin{figure}[t]
    \centering
    \begin{tikzpicture}[cls plot, scale=0.92]
    \foreach \n/\dx in {1/0,2/4.5,3/9} {
        \begin{scope}[shift={(\dx,0)}]
            \coordinate (O) at (0,0);
            \coordinate (u1) at (0,1);
            \coordinate (u2) at (1,0);
            \coordinate (u3) at (0,-1);
            \coordinate (u4) at (-1,\n);
            \CLSwall{u2}{u1}{clsgray1}
            \CLSwall{u1}{u4}{clsgray2}
            \foreach \i in {-1,0,1}
                \foreach \j in {-1,...,\n}
                    \node[lattice] at (\i,\j) {};
            \draw[chord] (u2) -- (u4);
            \ifnum\n=1
                \draw[sep] (-1.45,1.725)
                    node[lbl, anchor=south east] {$H$} -- (1.45,0.275);
            \fi
            \ifnum\n=3
                \draw[sep] (-1.1,2.65) -- (1.1,-0.65)
                    node[lbl, anchor=west] {$H$};
            \fi
            \foreach \v in {u2,u3,u4} {\CLSray{\v}}
            \CLSray[tested ray]{u1}
            \node[origin] at (O) {};
            \foreach \v in {u1,u2,u3,u4} {\node[gen] at (\v) {};}
            \node[lbl, anchor=west] at (0.15,1.1) {$u_{\rho_1}$};
            \node[lbl, anchor=west] at (1.35,0) {$u_{\rho_2}$};
            \node[lbl, anchor=east] at (-1.15,\n) {$u_{\rho_4}$};
            \node[lbl, anchor=west] at (0.15,-1.1) {$u_{\rho_3}$};
            \node[slbl, anchor=north] at (0,-1.65) {$n=\n$};
        \end{scope}
    }
\end{tikzpicture}

    \caption{The criterion of \Cref{prop:hyperplane-test} at the ray $\rho_1$ of $\Sigma_n$, for the marking of $-K$ by primitive generators. Only for $n = 2$, no line separating $u_{\rho_1}$ from the other vectors exists.}
    \label{fig:hirzebruch}
\end{figure}
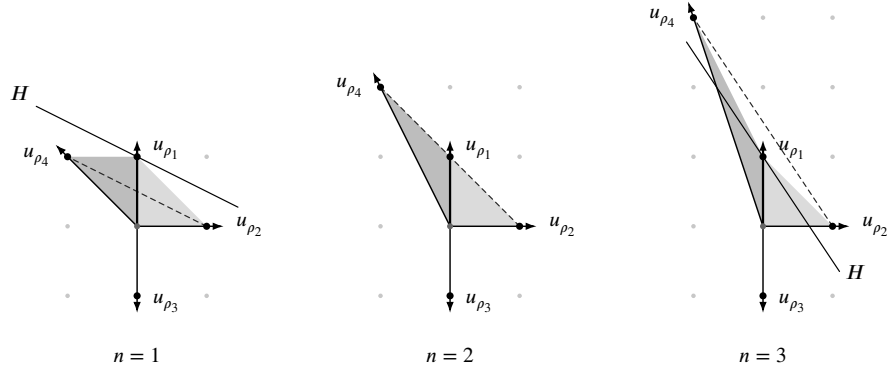

\begin{example}[Hirzebruch surfaces]\label{ex:hirzebruch}
    Let $F_n$ be the Hirzebruch surface associated to the complete fan $\Sigma_n$ in $\bb R^2$ with rays
    \[
    \rho_1 = \operatorname{cone}(0,1), \quad \rho_2 = \operatorname{cone}(1,0), \quad \rho_3 = \operatorname{cone}(0,-1), \quad \rho_4 = \operatorname{cone}(-1,n),
    \]
    so that $D_{\rho_2}, D_{\rho_4}$ are fibers of $F_n \to \bb P^1$, $D_{\rho_1}$ is the section of self-intersection $-n$, and $D_{\rho_3}$ the section of self-intersection $n$. Let us apply \Cref{thm:implementable} to the canonical divisor $K = -\sum_{i=1}^4 D_{\rho_i}$. By \Cref{rem:intermediate-cones}, only the rays require verification, and \Cref{fig:hirzebruch} depicts the equivalent criterion of \Cref{prop:hyperplane-test} at $\rho_1$ for $n = 1, 2, 3$. For $\tau = \rho_1$,
    \begin{align*}
        P_{K}(\rho_1)
        &= \{m \in M_{\bb R} \mid \bra{m, u_{\rho_1}} = 1, \ \bra{m, u_{\rho_2}} \geq 1, \ \bra{m, u_{\rho_4}} \geq 1\} \\
        &= \{(m_1, 1) \mid m_1 \geq 1, \ -m_1 + n \geq 1\} = [1, n-1] \times \{1\},\\
        P_{-K}(\rho_1)
        &= \{(m_1, -1) \mid m_1 \geq -1, \ -m_1 - n \geq -1\} = [-1, 1-n] \times \{-1\}.
    \end{align*}
    For $n \leq 1$, $P_K(\rho_1) = \emp$ and $\dim P_{-K}(\rho_1) = 1$. For $n \geq 3$, $\dim P_K(\rho_1) = 1$ and $P_{-K}(\rho_1) = \emp$. For $n = 2$, both polytopes are single points, so $\omega_{F_2}$ is not $\otimes$-generating: the obstruction is the $(-2)$-curve $V(\rho_1)$. Computations for the remaining rays give $P_{-K}(\rho_2) = \{-1\} \times [-1, 1]$, $P_{-K}(\rho_3) = [-1, n+1] \times \{1\}$, and $P_{-K}(\rho_4) = \{(nm_2 + 1, m_2) \mid -1 \leq m_2 \leq 1\}$, each of dimension $1$, with $P_K(\rho_i) = \emp$ for $i = 2,3,4$. Hence $\omega_{F_n}$ is $\otimes$-generating if and only if $n \neq 2$, recovering the case of Hirzebruch surfaces of \cite{ito2025derived}*{Corollary 5.14}.
\end{example}

\begin{example}[Toric Fano varieties]\label{ex:fano}
    Let us combinatorially show that the canonical bundle of a Gorenstein toric Fano variety is $\otimes$-generating. Recall that such a variety is $X = X_{\Sigma_P}$ for the normal fan $\Sigma_P$ of a reflexive polytope $P \subset M_{\bb R}$ \cite{cox2011toric}*{Theorem 8.3.4}, and that $\Sigma_P$ is the face fan of the polar dual $P^\circ \subset N_{\bb R}$, whose vertices are the primitive ray generators of $\Sigma_P$ \cite{cox2011toric}*{Exercise 2.3.4 and (8.3.2)}. The marking associated with $-K_X = \sum_\rho D_\rho$ is thus the set of vertices of $P^\circ$. A nonzero cone $\tau \in \Sigma_P$ is the cone over a proper face $F$ of $P^\circ$, with $\tau(1)$ the set of rays through the vertices of $F$, and the cones $\sigma \succ \tau$ are the cones over the faces containing $F$. Choose a supporting affine hyperplane $H$ of $P^\circ$ with $H \cap P^\circ = F$. Then $H$ contains the marking vectors of $\tau(1)$, and every other vertex of $P^\circ$, in particular every vertex of a face containing $F$ and not lying in $F$, lies strictly on the $P^\circ$-side of $H$. \Cref{prop:hyperplane-test} applies, and $\omega_X$ is $\otimes$-generating. By combinatorial means, this recovers the toric Gorenstein case of the statement that a variety with (ample or) anti-ample canonical bundle has $\otimes$-generating canonical bundle \cite{ito2025derived}*{Theorem 1.1}.
\end{example}

\subsection{\texorpdfstring{$\otimes$}{tensor}-generating fans: characterization of \texorpdfstring{$\otimes$}{tensor}-generating canonical bundles}\label{subsec:tensor-reflexive}

For a $\bb Q$-Gorenstein proper toric variety $X$ and an integer $l > 0$ with $lK_X$ Cartier, \Cref{cor:conv-aff} gives a criterion for $\otimes$-generation of $\ecal O_X(lK_X)$. In this section, we propose the definition of $\otimes$-generating fans suggested by this criterion, which generalizes the face fans of reflexive polytopes. 
In \Cref{sec:canonical}, we apply this criterion to blow-ups of $\bb P^3$ along configurations of invariant lines. 


\begin{definition}\label{def:tensor-reflexive}
    A complete fan $\Sigma$ in $N_{\bb R}$ is \emph{$\otimes$-generating} if for every cone $\tau \in \Sigma \setminus \{0\}$, some affine hyperplane in $N_{\bb R}$ contains the primitive generators $u_\rho$ with $\rho \in \tau(1)$ and has the primitive generators $u_\rho$, for the remaining rays of the cones containing $\tau$, strictly on one side.
\end{definition}

\begin{figure}[t]
    \centering
    \begin{tikzpicture}[cls plot, scale=1.05]
    \begin{scope}
        \coordinate (O) at (0,0);
        \coordinate (u0) at (-1,-1);
        \coordinate (u1) at (1,0);
        \coordinate (u2) at (0,1);
        \CLSwall{u1}{u2}{clsgray1}
        \CLSwall{u2}{u0}{clsgray2}
        \CLSwall{u0}{u1}{clsgray3}
        \foreach \i in {-1,0,1} \foreach \j in {-1,0,1}
            \node[lattice] at (\i,\j) {};
        \draw[slice] (u1) -- (u2) -- (u0) -- cycle;
        \foreach \v in {u0,u1,u2} {\CLSray{\v}}
        \node[origin] at (O) {};
        \foreach \v in {u0,u1,u2} {\node[gen] at (\v) {};}
        \node[lbl, anchor=south west] at (1.05,0.12) {$u_1$};
        \node[lbl, anchor=west] at (0.14,1.05) {$u_2$};
        \node[lbl, anchor=east] at (-1.14,-1.02) {$u_0$};
        \node[slbl, anchor=north] at (0,-1.65) {the fan of $\bb P^2$};
    \end{scope}
    \begin{scope}[shift={(4.4,0)}]
        \coordinate (O) at (0,0);
        \coordinate (u1) at (0,1);
        \coordinate (u2) at (1,0);
        \coordinate (u3) at (0,-1);
        \coordinate (u4) at (-1,3);
        \CLSwall{u2}{u1}{clsgray1}
        \CLSwall{u1}{u4}{clsgray2}
        \CLSwall{u4}{u3}{clsgray3}
        \CLSwall{u3}{u2}{clsgray2}
        \foreach \i in {-1,0,1} \foreach \j in {-1,0,1,2,3}
            \node[lattice] at (\i,\j) {};
        \draw[slice] (u2) -- (u1) -- (u4) -- (u3) -- cycle;
        \foreach \v in {u1,u2,u3,u4} {\CLSray{\v}}
        \node[origin] at (O) {};
        \foreach \v in {u1,u2,u3,u4} {\node[gen] at (\v) {};}
        \node[lbl, anchor=west] at (0.14,1.1) {$u_{\rho_1}$};
        \node[lbl, anchor=south west] at (1.05,0.12) {$u_{\rho_2}$};
        \node[lbl, anchor=west] at (0.16,-1.15) {$u_{\rho_3}$};
        \node[lbl, anchor=east] at (-1.14,3) {$u_{\rho_4}$};
        \node[slbl, anchor=north] at (0,-1.65) {the fan of $F_3$};
    \end{scope}
\end{tikzpicture}

    \caption{Two $\otimes$-generating complete fans in $\bb R^2$. The one on the left is reflexive, whereas the one on the right is not.}
    \label{fig:reflexive-fans}
\end{figure}
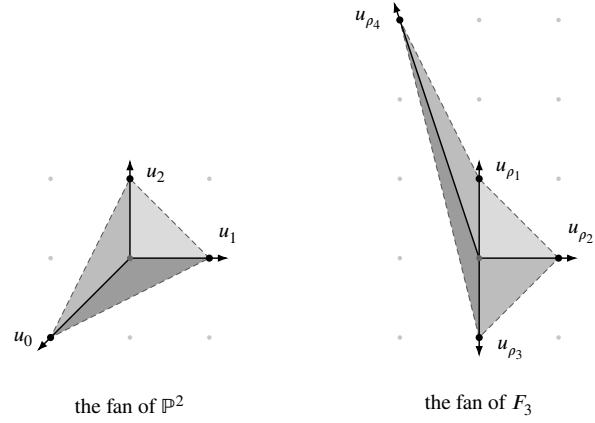

\Cref{fig:reflexive-fans} shows two $\otimes$-generating fans, with the left one reflexive and the right $\otimes$-generating but not reflexive. 

\begin{corollary}\label{cor:tensor-reflexive}
    Let $X = X_\Sigma$ be a proper toric variety. The following are equivalent.
    \begin{enumerate}
        \item The fan $\Sigma$ is $\otimes$-generating.
        \item The canonical divisor $K_X$ is $\bb Q$-Cartier, and $\ecal O_X(lK_X)$ is $\otimes$-generating for one, equivalently every, integer $l > 0$ with $lK_X$ Cartier.
    \end{enumerate}
    In particular, a Gorenstein proper toric variety has $\otimes$-generating canonical bundle if and only if its fan is $\otimes$-generating.
\end{corollary}
\begin{proof}
    We first show that the instances of \Cref{def:tensor-reflexive} at the maximal cones hold if and only if $K_X$ is $\bb Q$-Cartier. At a maximal cone $\sigma$, there are no remaining rays, so the condition asks only for an affine hyperplane through the generators $u_\rho$ with $\rho \in \sigma(1)$. Such a hyperplane cannot pass through the origin, since a linear hyperplane containing the generators of $\sigma$ would contain $\operatorname{span}(\sigma) = N_{\bb R}$. It can therefore be written $\{u \mid \bra{m_\sigma, u} = 1\}$ with $m_\sigma \in M_{\bb R}$, and since the system $\bra{m_\sigma, u_\rho} = 1$ for $\rho \in \sigma(1)$ has rational coefficients, its solvability is a rank condition, so a real solution exists if and only if a rational one does. Rational solutions for all maximal cones determine, after scaling and changing sign, Cartier data for a positive multiple of $-K_X = \sum_\rho D_\rho$, and conversely (\Cref{subsec:toric-conventions}).

    Now suppose $K_X$ is $\bb Q$-Cartier and let $l > 0$ be an integer with $lK_X$ Cartier. The divisor $D := -lK_X = \sum_\rho l D_\rho$ is Cartier with all coefficients positive, hence big by \Cref{cor:big-cone}, so \Cref{prop:hyperplane-test} applies to $\ecal O_X(-lK_X)$ with this $D$, whose marking is $u_{D, \rho} = u_\rho / l$. Rescaling by $l$ matches the hyperplanes for the marking with the hyperplanes for the primitive generators, so $\ecal O_X(-lK_X)$ is $\otimes$-generating if and only if \Cref{def:tensor-reflexive} holds at every nonzero cone that is not maximal. Since the maximal instances hold automatically when $K_X$ is $\bb Q$-Cartier, since a line bundle is $\otimes$-generating exactly when its inverse is, and since \Cref{def:tensor-reflexive} does not involve $l$, the equivalence of (i) and (ii) follows, for one $l$ and hence for every $l$. The Gorenstein case is $l = 1$.
\end{proof}

\begin{remark}\label{rem:reflexive-comparison}
Our definition of $\otimes$-generating fans generalizes the combinatorics of reflexive polytopes. When $X$ is a Gorenstein toric Fano variety, $\Sigma$ is the face fan of the reflexive polytope $P^\circ = \operatorname{conv}(u_\rho \mid \rho \in \Sigma(1))$, and a supporting affine hyperplane at the face corresponding to a nonzero cone $\tau$ contains the generators of $\tau(1)$ and has every other vertex strictly on one side, as in \Cref{ex:fano}. Face fans of reflexive polytopes are therefore $\otimes$-generating, and \Cref{cor:tensor-reflexive} recovers \Cref{ex:fano}.

   When $X$ is Gorenstein, the polytopes $\operatorname{conv}(u_\rho \mid \rho \in \sigma(1))$ for maximal cones $\sigma$ bound the star-shaped region
\[
B_\Sigma := \bigcup_\sigma \operatorname{conv}(\{0\} \cup \{u_\rho \mid \rho \in \sigma(1)\}),
\]
shown in \Cref{fig:reflexive-fans}. The variety $X$ is Fano exactly when $B_\Sigma$ is convex and these polytopes are its facets. For $\otimes$-generating fans, the separating hyperplanes in \Cref{def:tensor-reflexive} may have the origin on either side. Accordingly, $\omega_X$ may be big on some orbit closures and anti-big on others: on $F_n$ with $n \geq 3$, its restriction to the $(-n)$-section is big, while its restrictions to the other invariant curves are anti-big (\Cref{fig:hirzebruch}). In this language, \Cref{ex:hirzebruch} states that $\Sigma_n$ is $\otimes$-generating if and only if $n \neq 2$. Hence, the fan of a smooth projective toric surface is $\otimes$-generating precisely when there are no invariant $(-2)$-curves, agreeing with \cite{ito2025derived}*{Corollary 5.14}. We will see further examples in \Cref{thm:Xa,thm:Ya}, and \Cref{prop:P3-blowup-obstructions}.
\end{remark}

  Unlike reflexivity, $\otimes$-generation is a property of the fan rather than of the polytope spanned by the generators, and it is fragile: for example,  \Cref{thm:products} implies that $\otimes$-generation is not preserved under taking products. Furthermore, every nontrivial crepant modification destroys $\otimes$-generation, as follows.

\begin{corollary}\label{cor:crepant}
    Let $\Sigma' \neq \Sigma$ be a refinement of the complete fan $\Sigma$ and $\pi \colon X_{\Sigma'} \to X_\Sigma$ the induced toric morphism. Suppose the canonical divisors of $X_{\Sigma'}$ and $X_\Sigma$ are $\bb Q$-Cartier and $\pi$ is crepant, in the sense that $lK_{X_{\Sigma'}} \sim \pi^*(lK_{X_\Sigma})$ for an integer $l > 0$ with both sides Cartier. Then $\Sigma'$ is not $\otimes$-generating.
\end{corollary}
\begin{proof}
    We first produce a cone $\tau' \in \Sigma'$ of codimension one whose relative interior meets the interior of a maximal cone $\sigma \in \Sigma$. First, there exists some maximal cone $\sigma$ of $\Sigma$ that is not a cone of $\Sigma'$. Otherwise, each maximal cone of $\Sigma'$ would share interior points with some maximal cone of $\Sigma$, which would also be a cone of $\Sigma'$. Distinct cones of a fan have disjoint relative interiors, so the two cones would coincide, and the two fans would be equal. Every cone of $\Sigma'$ meeting $\operatorname{int}(\sigma)$ is contained in $\sigma$: it is contained in some cone $\tau$ of $\Sigma$, and the intersection $\tau \cap \sigma$ is a face of $\sigma$ containing an interior point, so $\sigma \subseteq \tau$, and $\tau = \sigma$ by maximality. The maximal cones of $\Sigma'$ contained in $\sigma$ therefore cover $\sigma$, and there are at least two of them, since $\sigma \notin \Sigma'$. Choose interior points $x$ and $y$ of two distinct such cones, so that the segment $[x, y]$ lies in the open convex set $\operatorname{int}(\sigma)$. For fixed $x$, the endpoints $y'$ for which $[x, y']$ meets a given cone of $\Sigma'$ of codimension at least $2$ form a subset of the cone over it with vertex $x$, so after replacing $y$ by a nearby interior point of its cone, the segment meets no such cone. Let $v$ be the point where $[x, y]$ leaves the maximal cone containing $x$. Then $v \in \operatorname{int}(\sigma)$, and $v$ lies in the relative interior of a codimension-one face $\tau'$ of that cone.

    Set $C := V(\tau') \subseteq X_{\Sigma'}$. The star fan of $\tau'$ is a complete fan in the line $N(\tau')_{\bb R}$, so $C \iso \bb P^1$. The set $\operatorname{relint}(\tau') \cap \operatorname{int}(\sigma)$ is nonempty, open in $\operatorname{relint}(\tau')$, and stable under positive scaling, so it contains a lattice point $w$. Since $w \in \operatorname{relint}(\tau')$, the limit $\lim_{t \to 0} \lambda^w(t)$ of the associated one-parameter subgroup exists in $X_{\Sigma'}$ and is the distinguished point of the orbit $O(\tau')$, and since $w \in \operatorname{int}(\sigma)$, the same limit computed in $X_\Sigma$ is the torus-fixed point $x_\sigma$ of the maximal cone $\sigma$. The morphism $\pi$ maps the distinguished point of $O(\tau')$ to $x_\sigma$, and by equivariance $\pi(O(\tau')) = \{x_\sigma\}$, so $\pi(C) = \{x_\sigma\}$.

    Crepancy gives $\ecal O_{X_{\Sigma'}}(lK_{X_{\Sigma'}}) \iso \pi^* \ecal O_{X_\Sigma}(lK_{X_\Sigma})$, and the restriction of the right-hand side to $C$ is pulled back along the constant map $C \to \{x_\sigma\}$, hence trivial of degree $0$. A line bundle on $\bb P^1$ is big or anti-big exactly when its degree is nonzero, so the restriction of $\ecal O_{X_{\Sigma'}}(lK_{X_{\Sigma'}})$ to the orbit closure $C$ is neither, the line bundle is not $\otimes$-generating by \Cref{thm:orbit-closure}, and $\Sigma'$ is not $\otimes$-generating by \Cref{cor:tensor-reflexive}.
\end{proof}

\begin{example}\label{ex:F2-P112}
    Let $P^\circ := \operatorname{conv}((1,0), (0,-1), (-1,2))$. The fan $\Sigma_2$ of $F_2$ from \Cref{ex:hirzebruch} refines the face fan of the triangle $P^\circ$ by the star subdivision at $u_{\rho_1} = (0,1)$, the primitive generator of $(1,0) + (-1,2)$. The triangle $P^\circ$ is reflexive, with facets on the lattice hyperplanes $\{x - y = 1\}$, $\{x + y = 1\}$, and $\{-3x - y = 1\}$, and its face fan is the fan of the weighted projective plane $\bb P(1,1,2)$, a Gorenstein toric Fano surface. The subdivision is crepant, since $(0,1)$ lies on the facet hyperplane $\{x + y = 1\}$ through $(1,0)$ and $(-1,2)$. \Cref{rem:reflexive-comparison} shows the fan of $\bb P(1,1,2)$ is $\otimes$-generating, \Cref{cor:crepant} shows $\Sigma_2$ is not, and the contraction realizes the obstruction: the $(-2)$-curve $V(\rho_1) \subseteq F_2$ maps to the singular point of $\bb P(1,1,2)$, and $\omega_{F_2}$ restricts to it with degree $0$.
\end{example}

Polar duality is an involution of the set of reflexive polytopes \cite{cox2011toric}*{\S 8.3}. We do not know whether a comparable structure organizes $\otimes$-generating fans.

\begin{question}\label{q:duality}
    Is there an involution of the class of $\otimes$-generating fans extending polar duality of reflexive polytopes?
\end{question}
\section{The curve criterion: an analogue of Kleiman's criterion}\label{sec:curves}

In this section, we prove a derived category analogue of Kleiman's criterion of ampleness in the case of toric varieties (\Cref{thm:B}): a line bundle $\mathcal{L}$ on a proper toric variety fails to be $\otimes$-generating if and only if there exists some irreducible complete curve $C$ for which $\deg \mathcal{L}|_C = 0$. We remark that in general, the curve cannot be taken torus-invariant (see \Cref{ex:diagonal}).  

\subsection{Proof of the curve criterion}

We use the following variant of Farkas' lemma. 
\begin{lemma}\label{lem:Farkas}
    Let $H \subset \bb R^n$ be a linear subspace and take $v \in \bb R^n$. Then, there exists $h \in H$ such that $v + h \in \bb R_{>0}^n$ if and only if $\lambda \cdot v > 0$ for any $0 \neq\lambda \in H^\perp \cap \bb R^n_{\geq 0}$.
\end{lemma}
\begin{proof}
    For the ``only if'' direction, suppose we have $h \in H$ with $v+h \in \bb R^n_{>0}$ and take $0 \neq\lambda \in H^\perp \cap \bb R^n_{\geq 0}$. Then, $0 < \lambda \cdot (v+h) = \lambda \cdot v$.

    Conversely, suppose $\lambda \cdot v > 0$ for any $0 \neq\lambda \in H^\perp \cap \bb R^n_{\geq 0}$ and set $A:= v + H$ to be the affine translation. Assume $A \cap \bb R^n_{>0} = \emp$ for contradiction. Since $A$ and $\bb R^n_{>0}$ are convex, there exist $0 \neq \lambda \in \bb R^n$ and $a \in \bb R$ such that $\lambda\cdot (v + h) \leq a \leq \lambda \cdot w$ for any $h \in H$ and any $w \in \bb R^n_{>0}$ by the hyperplane separation theorem. Note that we necessarily have $\lambda = (\lambda_i) \in \bb R_{\geq 0}^n$ since if $\lambda_i < 0$, then we can make $\lambda \cdot w$ arbitrarily negative. Hence, we can take $a = \inf_{w \in \bb R_{>0}^n} (\lambda \cdot w) = 0$. Moreover, $\lambda \cdot h = 0$ for any $h\in H$ since otherwise $\lambda\cdot t h$ is not bounded as $t \to \pm \infty$. Hence, we have $0 \neq\lambda \in H^\perp \cap \bb R^n_{\geq 0}$ with $\lambda \cdot v \leq 0$, which is a contradiction.
\end{proof}

\begin{prop}\label{prop:neither-big}
    Let $X = X_\Sigma$ be a proper toric variety and $D = \sum_\rho a_\rho D_\rho$ a $T_N$-invariant $\bb Q$-Cartier divisor. The following are equivalent.
    \begin{enumerate}
        \item Neither $D$ nor $-D$ is big.
        \item There is a nonzero vector $\lambda \in \bb R_{\geq 0}^{\Sigma(1)}$ with $\sum_\rho \lambda_\rho u_\rho = 0$ and $\sum_\rho \lambda_\rho a_\rho = 0$.
    \end{enumerate}
    The vector $\lambda$ in (ii) can be chosen with integer entries.
\end{prop}
\begin{proof}
    Set $a := (a_\rho)_\rho$ and
    \[
    H := \im\bigl(M_{\bb R} \ni m \mapsto (\bra{m,u_\rho})_\rho \in \bb R^{\Sigma(1)}\bigr),
    \qquad
    K := H^\perp \cap \bb R^{\Sigma(1)}_{\geq 0}.
    \]
    Then $K = \{\lambda \in \bb R^{\Sigma(1)}_{\geq 0} \mid \sum_\rho \lambda_\rho u_\rho = 0\}$. By \Cref{lem:toric-bigness,rem:Q-Cartier-bigness}, $D$ is big if and only if there exists $m \in M_{\bb R}$ such that $a_\rho + \bra{m,u_\rho} > 0$ for all $\rho \in \Sigma(1)$. By \Cref{lem:Farkas}, this is equivalent to $a \cdot \lambda > 0$ for every $0 \neq \lambda \in K$. Similarly, $-D$ is big if and only if $a \cdot \lambda < 0$ for every $0 \neq \lambda \in K$. Therefore, (ii) implies (i). 
    
    Conversely, if neither $D$ nor $-D$ is big, there are nonzero $\lambda^-,\lambda^+ \in K$ such that $a \cdot \lambda^- \leq 0 \leq a \cdot \lambda^+$. If either value is zero, we are done. Otherwise, a positive linear combination of $\lambda^-$ and $\lambda^+$ gives a nonzero $\lambda \in K$ with $a \cdot \lambda = 0$, proving (ii). Since the locus of $\lambda \in \bb R^{\Sigma(1)}_{\geq 0}$ satisfying the desired equalities is cut out by rational linear forms, it is a nontrivial rational polyhedral cone. Therefore, $\lambda$ can be chosen with integer entries. 
\end{proof}

\begin{prop}[Construction of the obstruction curve]\label{prop:degree-zero-curve}
    Let $X = X_\Sigma$ be a proper toric variety over $k$ and $D$ a $T_N$-invariant $\bb Q$-Cartier divisor such that neither $D$ nor $-D$ is big. Then there is an irreducible complete curve $C \subseteq X$ meeting the torus $T_N$ such that $\deg\bigl(\ecal O_X(lD)|_C\bigr) = 0$ for any integer $l > 0$ with $lD$ Cartier. Moreover, there is a finite morphism $f: \tilde C \to X$ from a normal irreducible projective curve $\tilde C$ such that $f^*\ecal O_X(lD)$ is torsion for any integer $l > 0$ with $lD$ Cartier.
\end{prop}
\begin{proof}
    By \Cref{prop:neither-big} we can take primitive ray generators $u_1,\ldots,u_s \in N$ and positive integers $c_1,\ldots,c_s$ such that
    \[
        \sum_{i=1}^s c_i u_i = 0
        \quad\text{and}\quad
        \sum_{i=1}^s c_i\varphi_D(u_i) = 0.
    \]
    Let $N' := \bb Z^s/\bb Z(1,\ldots,1)$, and let $e_1,\ldots,e_s$ be the images of the standard basis vectors. The first relation defines a lattice map $\psi:N'\to N$ by $\psi(e_i)=c_i u_i$. Let $\Sigma'$ be the fan with rays generated by the $e_i$ and no higher dimensional cones. Then $\psi$ is compatible with the fans $\Sigma'$ and $\Sigma$, hence inducing a toric morphism $\pi:X_{\Sigma'}\to X$.

    Let $l>0$ be any integer with $lD$ Cartier. The support function of $\pi^*(lD)$ is $l\varphi_D\circ\psi$ on $|\Sigma'|$ \cite{cox2011toric}*{Proposition 6.2.7}. The integers $lc_i\varphi_D(u_i)$ sum to zero, so they define $g\in M':=\operatorname{Hom}(N',\bb Z)$ with $g(e_i)=lc_i\varphi_D(u_i)$. Hence, the support function of $\pi^*(lD)$ is the restriction of $g$, so $\pi^*(lD)=\operatorname{div}(\chi^{-g})$ is principal and $\pi^*\ecal O_X(lD)\iso\ecal O_{X_{\Sigma'}}$.

    Since any $s-1$ rays form a strongly convex cone, we can form a complete fan $\Sigma''$ by adding cones spanned by fewer than $s$ rays, which gives a compactification $X_{\Sigma'}\subseteq X_{\Sigma''}$ with boundary of codimension at least $2$. Here $X_{\Sigma''}\iso\bb P^{s-1}$. Over an algebraic closure $\bar k$, if $s=2$, take $C':=X_{\Sigma',\bar k}\iso\bb P^1_{\bar k}$. If $s\geq3$, the intersection $C'$ of $s-2$ general hyperplanes in $\bb P^{s-1}_{\bar k}$ is a line avoiding $X_{\Sigma'',\bar k}\setminus X_{\Sigma',\bar k}$ and meeting each invariant divisor properly. In both cases, $C'$ is an irreducible complete curve in $X_{\Sigma',\bar k}$ meeting the torus and each invariant divisor. Since $\pi$ maps the torus into $T_N$ and the divisor associated to $e_i$ into the divisor associated to $u_i$, the image of $C'$ is not a point.

    Take a finite field extension $k\subset k'\subset\bar k$ so that $C'$ descends to a rational curve in $X_{\Sigma',k'}$, again denoted $C'$. Let $h:C'\to X$ be the composition with the toric map and projection to $X$. Its reduced image $C:=h(C')\subseteq X$ is an irreducible complete curve meeting $T_N$. The map $C'\to C$ is a finite surjective morphism of complete curves, so by the projection formula,
    \[
        \deg(C'/C)\cdot\deg_k\bigl(\ecal O_X(lD)|_C\bigr)
        = \deg_k\bigl(h^*\ecal O_X(lD)\bigr) = 0,
    \]
    and $\deg\bigl(\ecal O_X(lD)|_C\bigr) = 0$. Finally, $C'$ is a normal irreducible projective curve over $k$, so taking $\tilde C = C'$ and $f=h$ proves the last claim.
\end{proof}

\begin{theorem}\label{thm:curve-criterion}
    Let $X = X_\Sigma$ be a proper toric variety over $k$ and $\ecal L$ a line bundle on $X$. Then the following are equivalent:
    \begin{enumerate}
        \item $\ecal L$ is $\otimes$-generating,
        \item For every irreducible complete curve $C \subseteq X$, $\deg(\ecal L|_C) \neq 0$,
        \item $\ecal L$ is nowhere torsion in the sense of \cite{itoolander2026botamestack}, i.e., for every finite morphism $f: C \to X$ from a normal projective irreducible curve $C$, $f^*\ecal L$ is not a torsion line bundle.   \qedhere 
    \end{enumerate}
\end{theorem}
\begin{proof}
(i) $\implies$ (ii) follows by \cite{ito2025derived}*{Proposition 5.1}, noting that an irreducible complete curve in a proper variety is an integral projective curve. For (ii) $\implies$ (iii), note that the image of a normal projective irreducible curve $C$ under a finite morphism $f: C\to X$ is an irreducible complete curve $C' \subset X$, and if $\deg (\ecal L|_{C'}) \neq 0$, then $\deg f^*\ecal L \neq 0$ and hence is not torsion. Finally, we show the contrapositive of (iii) $\implies$ (i). Suppose $\ecal L$ is not $\otimes$-generating, and write $\ecal L \iso \ecal O_X(D)$ for a $T_N$-invariant Cartier divisor $D$. By \Cref{thm:implementable}, there is a cone $\tau \in \Sigma$ with $\dim P_D(\tau)$ and $\dim P_{-D}(\tau)$ both different from $\dim X - \dim\tau$. By \Cref{lem:restriction-polytope} and \Cref{lem:toric-bigness}, the restriction of $\ecal L$ to the proper toric variety $V(\tau)$ is $\ecal O_{V(\tau)}(D_\tau)$ with neither $D_\tau$ nor $-D_\tau$ big. Then, we conclude by applying \Cref{prop:degree-zero-curve} to $(V(\tau), D_\tau)$. 
\end{proof}

\begin{example}\label{ex:diagonal}
    The curve in \Cref{thm:curve-criterion} cannot in general be taken torus-invariant. Let $X = \bb P^1 \times \bb P^1$ and $\ecal L = \ecal O_X(1,-1)$. The restriction of $\ecal L$ to each of the four invariant curves is $\ecal O(1)$ or $\ecal O(-1)$, of nonzero degree, so no invariant subvariety of positive codimension witnesses the failure. Nevertheless, $\ecal L$ is neither big nor anti-big, since $\h^0(X, \ecal L^{\tens l}) = 0$ for all $l \neq 0$, so $\ecal L$ is not $\otimes$-generating. A curve as in \Cref{thm:curve-criterion} is the diagonal $\Delta \subset X$, of type $(1,1)$, with $\deg(\ecal L|_\Delta) = 1 - 1 = 0$. In the construction of \Cref{prop:degree-zero-curve}, the diagonal arises from the balanced pair $w_1 = (1,1)$, $w_2 = (-1,-1)$, on which the support function of the divisor of bidegree $(1,-1)$ restricts linearly.
\end{example}

\subsection{The movable cone of curves}\label{subsec:movable-cone}

In this subsection, we work over $k = \bb C$. For a complex projective variety $X$, we write $\overline{\operatorname{SME}}(X) \subseteq N_1(X)_{\bb R}$ for the closure of the cone of strongly movable curve classes, which by \cite{boucksom2013pseudo}*{Theorem 2.2} is dual to the pseudoeffective cone $\operatorname{PsEff}(X) \subseteq N^1(X)_{\bb R}$. 

\begin{corollary}\label{cor:sme-polyhedral}
    Let $X = X_\Sigma$ be a complex projective toric variety. Then $\overline{\operatorname{SME}}(X)$ is a rational polyhedral cone, and each of its extremal rays is spanned by the class of an irreducible curve meeting the torus.
\end{corollary}
\begin{proof}
    First, the cone of effective divisor classes in $N^1(X)_{\bb R}$ is spanned by the classes of the effective $T_N$-invariant Cartier divisors: a nonzero space of sections of a line bundle on $X$ decomposes into character eigenspaces, and the divisor of an eigensection is invariant. The invariant effective Cartier divisors form the intersection of the rational polyhedral cone $\bb R^{\Sigma(1)}_{\geq 0}$ with the rational subspace of Cartier $\bb R$-divisors, so the effective cone is the image of a rational polyhedral cone under a linear map defined over $\bb Q$. It is therefore rational polyhedral and closed, and equals $\operatorname{PsEff}(X)$. Hence $\overline{\operatorname{SME}}(X) = \operatorname{PsEff}(X)^\vee$ is rational polyhedral as well.

    The cone $\operatorname{PsEff}(X)$ is full-dimensional, since it contains the open cone of ample classes. It is also strongly convex. To see this, suppose $\alpha$ and $-\alpha$ are both pseudoeffective, and write them as classes of effective invariant $\bb R$-Cartier divisors $E$ and $E'$ as in the previous paragraph, so that $E + E' \equiv 0$. Intersecting with $H^{n-1}$ for an ample divisor $H$ gives $\deg_H E + \deg_H E' = 0$ with both degrees nonnegative, so $\deg_H E = 0$. Since every prime invariant divisor has positive $H$-degree, every coefficient of $E$ vanishes and $\alpha = 0$.

    Let $R$ be an extremal ray of $\overline{\operatorname{SME}}(X) = \operatorname{PsEff}(X)^\vee$ and $F \subset \operatorname{PsEff}(X)$ the corresponding facet, so that $R = \operatorname{PsEff}(X)^\vee \cap F^\perp$. Here a class $\gamma \in \operatorname{PsEff}(X)^\vee$ that vanishes on a relative interior point of $F$ vanishes on all of $F$, since $F$ is spanned by finitely many generators appearing with positive coefficients in that point and $\gamma$ is nonnegative on each. Choose a rational class in the relative interior of $F$, represented by a $T_N$-invariant $\bb Q$-Cartier divisor $D$. The class of $D$ lies on the boundary of $\operatorname{PsEff}(X)$, so $D$ is not big, and $-D$ is not big either: a big $-D$ would place $-[D]$ in the interior of $\operatorname{PsEff}(X)$, contradicting strong convexity. \Cref{prop:degree-zero-curve} now gives an irreducible curve $C$ meeting $T_N$ with $D \cdot C = 0$. Since $C$ meets the torus, $C$ is contained in no invariant divisor, so $C$ pairs nonnegatively with every effective invariant $\bb R$-Cartier divisor, and therefore with $\operatorname{PsEff}(X)$. Thus $[C] \in \operatorname{PsEff}(X)^\vee \cap D^\perp = R$, and $[C] \neq 0$ because $C \cdot H > 0$. Hence $R$ is spanned by the class of $C$.
\end{proof}

\section{Existence and non-existence of \texorpdfstring{$\otimes$}{tensor}-generating line bundles}\label{sec:existence}

Every projective variety has a $\otimes$-generating line bundle, namely any ample one. It is well known that there are smooth proper toric varieties that are not projective \cite{oda1978torus}*{Theorem 9.6}, so it is natural to ask whether every proper toric variety still admits a $\otimes$-generating line bundle. In this section, we show that the answer is positive for simplicial varieties of dimension at most three and negative in general: products of complete toric varieties admit a $\otimes$-generating line bundle only when they are projective. This produces the first smooth proper varieties known to us with no $\otimes$-generating line bundle.

\subsection{Existence in dimension at most three}\label{subsec:existence-low-dim}

\begin{prop}\label{prop:existence-threefolds}
    Let $X = X_\Sigma$ be a simplicial proper toric variety of dimension $n \leq 3$. Then $X$ admits a $\otimes$-generating line bundle.
\end{prop}
\begin{proof}
    Fix an inner product on the underlying lattice of $\Sigma$. Let
    \[
D_0=\sum_{\rho\in\Sigma(1)}\|u_\rho\|D_\rho
    \]
    be an $\bb R$-divisor whose marking vectors $\frac{u_\rho}{\|u_\rho\|}$ lie on the unit sphere. Our goal is to show that a small perturbation of $D_0$ would pass the test in \Cref{prop:hyperplane-test} or \Cref{cor:conv-aff}. Since we can take the tangent hyperplane of the unit sphere at a marking vector to separate it from other vectors, it remains to run the test on $2$-dimensional cones in a $3$-dimensional fan. The strict separations at the rays persist under sufficiently small perturbations.

    Since a codimension $1$ cone is strictly contained in only two other cones, we are looking at the intersection of the convex hull of two points and a line. Write $u_1,u_2$ for the primitive generators of such a cone and $u_3,u_4$ for the other generators of its two adjacent maximal cones. These vectors have rank $3$, so their relations are generated by a nonzero relation $\sum_{i=1}^4 c_i u_i=0$. If $a_i>0$ are the perturbed coefficients of these rays, the only linear relations among $u_i/a_i$ are proportional to $(c_i a_i)_i$. Hence these four points are affinely dependent precisely when
    \[
        \sum_{i=1}^4 c_i a_i=0.
    \]
    If the line through the first two points meets the segment joining the other two, this equation holds.

    There are only finitely many such proper hyperplanes in the coefficient space, whose complement is dense and open. Therefore, we can perturb $D_0$ to make it a $\bb Q$-divisor that passes the test for all cones. This divisor is $\bb Q$-Cartier since $\Sigma$ is simplicial, and big since its coefficients are positive. Let $D$ be a positive Cartier multiple of this perturbation, and $\ecal O_X(D)$ would be a $\otimes$-generating line bundle.
\end{proof}

\subsection{\texorpdfstring{$\otimes$}{tensor}-generation on products}\label{subsec:products}

We now prove a somewhat surprising fact about the existence of $\otimes$-generating line bundles on products of proper varieties. That is, if $Y := \prod_{i=1}^n X_i$ is a product of positive-dimensional proper geometrically integral varieties with $n \geq 2$, then $Y$ admits a $\otimes$-generating line bundle if and only if each $X_i$ is projective. Furthermore, as long as the Picard group of $Y$ splits as a direct sum of $\operatorname{Pic}(X_i)$, we can conclude that $\otimes$-generation on $Y$ is equivalent to ampleness or anti-ampleness. For example, $Y$ will have this property if factors are normal, geometrically integral, and rational, which includes all toric varieties with our convention. Along these lines, we also give a combinatorial proof for the case of products of complete toric varieties in \Cref{thm:product-criterion-combinatorial}. We first provide some lemmas on bigness.

\begin{lemma}\label{lemma:bignessdescend}
    Let $X$ be a proper geometrically integral variety over a field $k$. For any field extension $k \subset K$, a line bundle $\ecal L$ on $X$ is big if and only if so is the pullback $\ecal L_K$ of $\ecal L$ to $X_K$.
\end{lemma}
\begin{proof}
    By the flat base change, we have
    \[
    \h^0(X, \ecal L^{\otimes n})  \otimes_k K = \h^0(X_K, \ecal L_K^{\otimes n})
    \]
    for every $n \geq 0$ and therefore $\dim_K \h^0(X_K, \ecal L_K^{\otimes n}) = \dim _k \h^0(X, \ecal L^{\otimes n})$. Since $X_K$ is proper and integral with $\dim X_K =\dim X$ by assumption, we conclude by definition (i.e., characterization (vii) of \cite{ito2025derived}*{Lemma 2.9}).
\end{proof}

\begin{lemma}\label{lem:invariance-of-bigness}
    Let $X$ and $Y$ be proper geometrically integral varieties over a field $k$ and let $\ecal L$ be a big line bundle on $X \times Y$. Then, for any point $x \in X$, $\ecal L|_{\{x\} \times Y}$ is big.
\end{lemma}

\begin{proof}
Let $\bar k$ be an algebraic closure of $k$. We first show that the claim holds for all $k$-points in $X$. By \Cref{lemma:bignessdescend} and our assumptions, we may assume $k = \bar k$. First, we claim that for any closed points $x_1, x_2 \in X$, $\ecal L|_{\{x_1\}\times Y}$ and $\ecal L|_{\{x_2\}\times Y}$ are numerically equivalent under the canonical isomorphism $\{x_1\} \times Y \cong Y \cong \{x_2\} \times Y$. Indeed, for any integral closed curve $C \subset Y$, we have 
\[
\deg_C \ecal L|_{\{x_1\}\times C} = \chi(C, \ecal L|_{\{x_1\}\times C}) - \chi(C, \ecal O_C) = \chi(C, \ecal L|_{\{x_2\}\times C}) - \chi(C, \ecal O_C) = \deg_C \ecal L|_{\{x_2\}\times C}
\]
by the local constancy of the Euler characteristic along $X \times C \to X$ (\cite{stacks-project}*{\href{https://stacks.math.columbia.edu/tag/0B9T}{Tag 0B9T}}). Since bigness is a numerical property over an algebraically closed field by \cite{ito2025derived}*{Lemma 4.9} (see also Lemma 4.7 in loc. cit.), it suffices to show that there exists a closed point $x \in X$ such that $\ecal L|_{\{x\} \times Y}$ is big. Now, since $\ecal L$ is big, there exists $s \in \h^0(X\times Y, \ecal L^{\otimes n})$ for some $n \geq 0$ such that $(X\times Y)_s \subset X\times Y$ is a non-empty and affine open subscheme by \cite{ito2025derived}*{Lemma 2.8}. Take a closed point $x \in X$ in the image of the projection $(X\times Y)_s \to X$. Then, for $s|_{\{x\}\times Y} \in \h^0(\{x\}\times Y, (\ecal L|_{\{x\} \times Y})^{\otimes n})$, we have that 
\[
(\{x\}\times Y)_{s|_{\{x\}\times Y}} = (X\times Y)_s \cap (\{x\} \times Y) 
\]
is a non-empty and affine open subscheme of $\{x\} \times Y$ and therefore $\ecal L|_{\{x\} \times Y}$ is big again by \cite{ito2025derived}*{Lemma 2.8}. Now, we assume $k$ is an arbitrary field again and take any point $x \in X$. Set $K = \kappa(x)$ and take a canonical lift $\tilde x \in X_K(K)$. Since $\ecal L_K$ is big on $(X\times Y)_K \cong X_K \times_K Y_K$ by \Cref{lemma:bignessdescend} and $X_K$ and $Y_K$ are proper geometrically integral varieties over $K$, the first part of the proof shows that $\ecal L_K|_{\{\tilde x\} \times Y_K}$ is big. Since $\ecal L_K|_{\{\tilde x\} \times Y_K} \cong \ecal L|_{\{x\} \times Y}$ under the canonical isomorphism $\{\tilde x\} \times_K Y_K \cong \{x\} \times Y$, we are done.
\end{proof}

\begin{theorem}\label{thm:products}
Let $n \geq 2$ and let $X_1, \dots, X_n$ be positive-dimensional proper geometrically integral varieties. The following are equivalent for $Y := \prod_{i=1}^n X_i$.
    \begin{enumerate}
        \item $Y$ admits a $\otimes$-generating line bundle.
        \item Every $X_i$ is projective.
    \end{enumerate}
    Moreover, a line bundle of the form $\ecal L_1 \boxtimes \cdots \boxtimes \ecal L_n$ on $Y$ for $\ecal L_i \in \pic(X_i)$ is $\otimes$-generating if and only if it is ample or anti-ample. In particular,
   if {the external tensor product map $\bigoplus_i \pic(X_i) \to \pic(Y)$ is furthermore an isomorphism (e.g. when $X_i$'s are toric by \cite{gonzalezaviles2017units}*{Corollary 5.7} )}, then a line bundle on $Y$ is $\otimes$-generating if and only if it is ample or anti-ample.
\end{theorem}


\begin{proof}
    It suffices to consider the case when $n=2$. Furthermore, we may assume that $k$ is algebraically closed by \cite{ito2025derived}*{Lemma 2.5} and \cite{stacks-project}*{\href{https://stacks.math.columbia.edu/tag/0BDG}{Tag 0BDG}}.  
    
    (ii) $\implies$ (i). If every $X_i$ is projective, then $Y$ is projective, and an ample line bundle is $\otimes$-generating by \cite{ito2025derived}*{Theorem 1.1}, since its restriction to every closed subvariety is again ample, hence big.

 (i) $\implies$ (ii). Now assume $\ecal L$ is $\otimes$-generating. Then, the restriction of $\ecal{L}$ is either big or anti-big on every closed subvariety of $Y$. We can without loss of generality assume that $\ecal{L}$ is big on $Y$. 
 Then, by \Cref{lem:invariance-of-bigness}, for all $x\in X_1$, $\ecal{L}|_{x\times X_2}$ is big on $x\times X_2$.  Suppose for contradiction that there exists a subvariety $Z\subseteq X_1$ for which $\ecal{L}|_{Z\times X_2}$ is anti-big. Again by \Cref{lem:invariance-of-bigness}, $\ecal{L}|_{x\times X_2}$ must be anti-big for all points $x\in Z \subseteq X_1$,  which is a contradiction,
 {since then, after raising $\ecal L$ to a suitable power, $\ecal{L}|_{x\times X_2}$ is isomorphic to the structure sheaf as both $\ecal{L}|_{x\times X_2}$ and $\ecal{L}\inv|_{x\times X_2}$ have non-zero sections, whose product must be a non-zero constant function. But $X_2$ is positive dimensional so the structure sheaf is not big by definition}. Therefore, $\ecal{L}|_{Z\times X_2}$ is big for all subvarieties $Z\subseteq X_1$. By applying \Cref{lem:invariance-of-bigness} to the other factor, we see that for a fixed copy of $X_1 = X_1\times x_2$ and for all subvarieties $Z = Z\times x_2$, the restriction $\ecal{L}|_{Z}$ is big on $Z$. Nakai--Moishezon then implies that $\ecal{L}|_{X_1}$ is ample. The same proof goes for $X_2$. 

    For the last claim, the argument just given shows that a $\otimes$-generating $\ecal L$ has, after possibly passing to $\ecal L^{-1}$, all $\ecal L_i$ ample, so $\ecal L = \ecal L_1 \boxtimes \cdots \boxtimes \ecal L_n$ is ample. Conversely, ample and anti-ample line bundles are $\otimes$-generating as in (ii) $\implies$ (i).

    For the final claim, assume that the Picard group of $Y$ factors as
    \[
        \mathrm{Pic}(Y)=\bigoplus_i \pi_i^*\mathrm{Pic(X_i)}.
    \]
    Then, every line bundle $\ecal L$ on $Y$ can be written as the box product of line bundles on each factor, which can be determined by restricting $\ecal L$ to any fiber. When $\ecal L$ is $\otimes$-generating and big, the argument above shows that these restrictions are ample, and therefore $\ecal L$ is itself ample.
\end{proof}
\begin{remark}
    For complete toric varieties, we have the following simpler arguments. Let $X_1,X_2$ be complete positive-dimensional toric varieties over a field. We can see that if a line bundle $\ecal L$ on $X_1 \times X_2$ is $\otimes$-generating, then it is either ample or anti-ample. As above, we may work over an algebraically closed field and we may write $\ecal L \cong \ecal L_1 \boxtimes \ecal L_2$ for line bundles $\ecal L_1$ on $X_1$ and $\ecal L_2$ on $X_2$. Let $C_1 \subseteq X_1$ and $C_2 \subseteq X_2$ be invariant curves, which exist since each factor is complete of positive dimension. Set
    \[
    S := C_1 \times C_2  \subset X_1 \times X_2,
    \]
  which is isomorphic to $\bb P^1 \times \bb P^1$. The restriction $\ecal L|_S$ has bidegree $(\deg_{C_1}\ecal L_1|_{C_1}, \deg_{C_2}\ecal L_2|_{C_2})$. Since ampleness and bigness agree on $\bb P^1 \times \bb P^1$, a line bundle of bidegree $(p, q)$ on $\bb P^1 \times \bb P^1$ is big if and only if $p > 0$ and $q > 0$ and anti-big if and only if $p < 0$ and $q < 0$. Since $\ecal L$ is $\otimes$-generating and $S$ is invariant, $\ecal L|_S$ is big or anti-big by \Cref{thm:orbit-closure} (or \cite{ito2025derived}*{Theorem 1.1}), so
    \[
    (\deg_{C_1}\ecal L_1|_{C_1}) (\deg_{C_2}\ecal L_2|_{C_2}) > 0.
    \]
    In particular, the sign of $\deg_C \ecal L_1|_C$ is consistently positive or negative for all invariant curves $C \subset X_1$. By the toric Kleiman criterion \cite{cox2011toric}*{Theorem 6.3.13}, $\ecal L_1$ is either ample or anti-ample. Similarly, $\ecal L_2$ is either ample or anti-ample and by the inequality above, $\ecal L_1$ and $\ecal L_2$ are both ample or both anti-ample.  
\end{remark}
 We remark that the above proof can fail if we consider fiber bundles instead of products. See \Cref{ex:bundle-over-F4}.  The facts above can also be deduced from the combinatorics of marked fans.  First, we need to slightly generalize the definition of tensor-generating fans (\Cref{def:tensor-reflexive}).

\begin{definition}\label{def:tg-marking}
    A complete marked fan $(\Sigma,u )$ inside a real vector space $V$ (no integral structure required) is \textit{$\otimes$-generating} if for every cone $\tau\in \Sigma \setminus \{0\}$, some affine hyperplane in $V$ contains the marking vectors $u_{\rho}$ when $\rho\in\tau(1)$, and has the marking vectors $u_{\rho}$, for the remaining rays of the cones containing $\tau$, strictly on one side. If the affine hyperplanes can be chosen such that the marking vectors for the remaining rays and the origin stay strictly on one side, the marking is said to be \textit{ample}.
\end{definition}

When $\Sigma$ is rational, the anti-canonical marking is $\otimes$-generating if and only if $\Sigma$ is $\otimes$-generating. When the marking is the associated marking $u_D$ of a big $\bb Q$-Cartier divisor $D=\sum_\rho a_\rho D_\rho$, we've seen in \Cref{prop:hyperplane-test} that the marking is $\otimes$-generating if and only if $\mathscr{O}(lD)$ is $\otimes$-generating for some integer $l$. The relation between ample markings and ample divisors is depicted by the following lemma.

\begin{lemma}\label{lem:ample-markings}
    A marking $u$ on a fan $\Sigma$ is ample if and only if they form the vertices of a convex polytope, whose face fan is precisely $\Sigma$.
\end{lemma}
\begin{proof}
    We only need to consider $\dim V>1$ since the case when $\dim V=1$ is trivial. For the if direction, let $\Sigma$ be the face fan of the polytope $P$ whose vertices are $\{u_\rho\}_{\rho\in\Sigma(1)}$. Then for each cone $\tau$ of $\Sigma$, the corresponding face on $P$ admits a supporting hyperplane. By definition, this hyperplane contains all the marking vectors for rays in $\tau$, and has all other marking vectors, along with the origin, on one side.

    For the only if direction, we first assume that the affine plane supporting cone $\tau$ is given by the linear equation $f_\tau=-1$, and hence the distinguished side having the origin is given by $f_\tau>-1$. We readily see that these $f_\tau|_\tau$'s glue together to form a piecewise linear support function $\varphi_u$ of $\Sigma$. We can hence define
    \[
        B=\bigcup_{\{0\}\neq\tau\in\Sigma}(\tau\cap\{f_\tau\geq-1\})=\{\varphi_u\geq-1\},
    \]
    which we claim to be a convex polytope with face fan $\Sigma$. Equivalently, we claim that $\varphi_u$ is strictly convex, since then, the intersection $B\cap \{f_\tau=-1\}$ would be equal to the convex hull of $u_\rho$ with $\rho\in \tau(1)$, implying that $\tau$ is in the face fan of $B$.
    
    According to the implication from $(c)$ to $(a)$ in \cite{cox2011toric}*{Lemma 6.1.13}, whose proof does not require an underlying integral structure of the fan, it suffices to show that for every wall $\tau=\sigma\cap\sigma'$, there is $\rho_0\in\sigma'\setminus\sigma$ where $\varphi_u(\rho_0)<f_\sigma(\rho_0)$. Note that $f_\sigma,f_\tau$ and $f_{\sigma'}$ agree on $\mathrm{span}(\tau)$. We pick arbitrary marking vectors $\rho_0\in \sigma'\setminus\sigma$ and $\rho_1\in\sigma\setminus\sigma'$. Since $V/\mathrm{span}(\tau)$ is $1$-dimensional, the images $\bar\rho_0,\bar\rho_1$ are negative multiples of each other. According to the assumption, we know that
    \[
        0<f_\tau(\rho_1)-f_\sigma(\rho_1)=(f_\tau-f_\sigma)(\bar\rho_1),
    \]
    and therefore,
    \[
    0>(f_\tau-f_{\sigma})(\bar\rho_0)=f_\tau(\rho_0)-f_\sigma(\rho_0),
    \]
    which implies
    \[
    \varphi_u(\rho_0)=f_{\sigma'}(\rho_0)<f_\tau(\rho_0)<f_\sigma(\rho_0).
    \]
\end{proof}

When the marking is the associated marking $u_D$ of a big $\bb Q$-Cartier divisor $D=\sum_\rho a_\rho D_\rho$, the support function we constructed is precisely the support function $\varphi_D$ of $D$. Therefore,  according to \cite{cox2011toric}*{Theorem 6.1.14}, \Cref{lem:ample-markings} shows that the divisor $D$ is ample if and only if the marking $u_D$ is ample.

Given two marked fans $(\Sigma,u),\ (\Sigma',u')$, since the rays of $\Sigma\times\Sigma'$ are in natural bijection with $\Sigma(1)\cup \Sigma'(1)$, we can form a marking by joining $u$ and $u'$, denoted as $u\oplus u'$. Conversely, every marking on the product fan can be written as a direct sum. We now prove the combinatorial analogue of \Cref{thm:products}.

\begin{theorem} \label{thm:product-criterion-combinatorial}
    Let $n>1$ be an integer and $(\Sigma_1,u_1),\ldots,(\Sigma_n,u_n)$ be marked fans of positive dimensions. Then, TFAE:
    \begin{enumerate}
        \item the direct sum $\bigoplus_i u_i$ on $\prod_i \Sigma_i$ is $\otimes$-generating;
        \item every marking $u_i$ is ample;
        \item the direct sum $\bigoplus_i u_i$ on $\prod_i \Sigma_i$ is ample. \qedhere 
    \end{enumerate} 
\end{theorem}

\begin{proof}
It suffices to prove the theorem when $n=2$. For simplicity, we denote the underlying vector spaces by $V_1,V_2$ and regard them as subspaces of $V_1\oplus V_2$.

We use \Cref{lem:ample-markings} to rephrase the equivalence of (ii) and (iii). When $\Sigma_i$ is the face fan of $P_i$, we can construct its \textit{direct sum} $P_1\oplus P_2=\mathrm{conv}(P_1\times0\cup0\times P_2)$, whose proper faces are of the form $\mathrm{conv}(F_1\times0\cup0\times F_2)$ for proper faces $F_i$ of $P_i$, allowing one empty factor. It's clear that the face fan of $P_1\oplus P_2$ is $\Sigma_1\times\Sigma_2$, and the vertices of $P_1\oplus P_2$ is the disjoint union of the vertices in $P_1$ and $P_2$; when $\Sigma_1\times\Sigma_2$ is the face fan of a polytope $P$, since $\Sigma_i$ is a subfan of $\Sigma_1\times\Sigma_2$ with support $V_i$, the face fan of the convex polytope $P\cap V_i$ must be $\Sigma_i$, and the marking given by the vertices in $P\cap V_i$ is ample.
    
(iii) implies (i) trivially, and we're left to show that (i) implies (ii). Let $\tau\in \Sigma_1\setminus\{0\}$ and $H$ be a hyperplane satisfying the conditions in \Cref{def:tg-marking}. Since for any $\rho\in\Sigma_2(1)$, we have $\rho\in(\tau\times\rho)(1)\setminus \tau(1)$, the distinguished side of $H$ contains all $u_{2,\rho}$. Since $\Sigma_2$ is complete of positive dimension, the convex hull of a set of ray generators must contain the origin, and therefore the origin also lies on the distinguished side of $H$. This shows that $u_1$ is ample.
\end{proof}

\begin{remark}
    The proof of (i) implying (ii) is completely analogous to that of \Cref{thm:products}. Looking at $\tau\in\Sigma_1\setminus\{0\}$ is analogous to looking at a subvariety of the first factor of a product, and the convex hull argument should be thought of as applying \Cref{lem:invariance-of-bigness} to the second factor.\end{remark}

\begin{corollary}\label{cor:no-TG}
    Let $X$ be a complete toric variety that is not projective. Then $X \times \bb P^1$ admits no $\otimes$-generating line bundle. In particular, if $X$ is a smooth non-projective complete toric threefold, then $X \times \bb P^1$ is a smooth proper toric fourfold with no $\otimes$-generating line bundle.
\end{corollary}
\begin{proof}
    Apply \Cref{thm:products} to $Y = X \times \bb P^1$.
\end{proof}

Smooth non-projective complete toric threefolds are classical: there are examples in the classification of \cite{oda1978torus}*{Theorem 9.6}, and the family $Z''_{13}$ studied in \Cref{sec:canonical} provides an infinite family of such. By \Cref{prop:existence-threefolds}, a smooth example of a proper toric variety with no $\otimes$-generating line bundle must have dimension at least $4$, so the fourfolds of \Cref{cor:no-TG} have the smallest possible dimension.

Without smoothness, examples with no $\otimes$-generating line bundle occur already in dimension $3$ and Picard number $1$, as the following observation shows.

\begin{example}[Fujino's example]\label{ex:fujino}
    Fujino \cite{fujino2005kleiman}*{Lemma 2.2} constructs a complete non-projective toric threefold $X_P$ with Picard number $1$ and canonical Gorenstein singularities, together with a small projective toric partner: a toric birational morphism $f \colon X_P \to X_Q$ to a projective toric threefold contracting a single curve $C \iso \bb P^1$. For an ample divisor $H$ on $X_Q$, the pullback $f^* H$ is a nontrivial nef Cartier divisor with $f^*H \cdot C = 0$. Since the Picard number of $X_P$ is $1$, every line bundle on $X_P$ is numerically proportional to $\ecal O_{X_P}(f^* H)$, and therefore has degree $0$ on $C$. By \Cref{thm:curve-criterion}, or directly by \cite{ito2025derived}*{Theorem 1.1}, no line bundle on $X_P$ is $\otimes$-generating.
\end{example}

    For smooth $X$, it is natural to consider the minimal dimension and Picard number such that it admits no $\otimes$-generating line bundle.   Smooth complete toric varieties of Picard number at most $3$ are projective in every dimension \cite{kleinschmidt1991smooth}. In dimension $3$, the unique non-projective smooth complete toric threefold of Picard number $4$ is the variety labeled $[7\text{-}5]$ in the classification \cite{oda1978torus}*{Theorem 9.6}, as observed by Fujino and Sato \cite{fujino2024non}*{Remark 1.6 and Example 4.1}. Taking $X = [7\text{-}5]$ in \Cref{cor:no-TG} therefore produces a smooth proper toric fourfold of Picard number $5$ with no $\otimes$-generating line bundle. The minimal dimension is therefore $4$ and the minimal Picard number is either $4$ or $5$. Candidates of Picard number $4$ could be sought in the classification of Choi--Jang--Vall\'ee \cite{choi2025complete}.

\begin{question}\label{q:minimal-picard}
    Is there a smooth proper toric variety with Picard number $4$ with no $\otimes$-generating line bundle?
\end{question}

\section{Examples of \texorpdfstring{$\otimes$}{tensor}-generating canonical bundles}\label{sec:canonical}

Following \cites{ito2025polarizations, ito2025derived}, a motivation to understand varieties with $\otimes$-generating canonical bundles is the Bondal--Orlov program: 
$\otimes$-generation of $\omega_X$ rules out non-isomorphic Fourier--Mukai partners and forces every Fourier--Mukai autoequivalence of $\perf(X)$ to be standard. We use the derived Nakai--Moishezon criterion (\Cref{thm:implementable}) in \S6.1, the $\otimes$-generating fan criterion (\Cref{cor:tensor-reflexive}) in \S6.2, and the derived Kleiman criterion (\Cref{thm:curve-criterion}) together with the rigidity on products (\Cref{thm:products}) in \S6.3, respectively. 



\subsection{The family \texorpdfstring{$Z''_{13}$}{Z''13}}\label{subsec:Z13}
In this subsection, we describe examples of smooth non-projective proper toric threefolds with $\otimes$-generating canonical bundles, using the combinatorial criterion of \Cref{thm:implementable}. To the authors' knowledge, these are first examples of (smooth) non-projective proper varieties without any projective Fourier--Mukai partner. They give an affirmative answer to \cite{ito2025polarizations}*{Question 4.8}.
\begin{construction}\label{const:Z13}
    For integers $a, b, c, d$, let $Z''_{13}(a,b,c,d)$ be the toric threefold whose fan $\Sigma \subset \bb R^3$ has rays generated by
    \begin{align*}
    v_1&=(1,1,b), & v_2&=(1,0,0), & v_3&=(0,-1,0), & v_4&=(0,0,-1),\\
    v_5&=(-1,a,d), & v_6&=(0,1,0), & v_7&=(0,0,1), & v_8&=(-1,c,d+1),
    \end{align*}
    and whose maximal cones are
    \[
    124,\ 234,\ 345,\ 456,\ 146,\ 167,\ 127,\ 237,\ 358,\ 568,\ 378,\ 678,
    \]
    where $ijk$ denotes $\operatorname{cone}(v_i, v_j, v_k)$.
\end{construction}

\begin{remark}\label{rem:Z13-attribution}
    The variety $Z''_{13}(a,b,c,d)$ is a smooth complete toric threefold, and the notation follows the classification of smooth complete toric threefolds of Picard number at most $5$ by Oda and Miyake \cite{oda1978torus}*{Theorem 9.6}, in the case labeled $[8\text{-}13'']$ there. The family is studied by Fujino and Sato \cite{fujino2024non}*{\S 4}, who show that $Z''_{13}(a,b,c,d)$ is projective if $b = 0$ or $a = c$, and non-projective if $b \neq 0$ and $a \neq c$.
\end{remark}
See the \hyperref[par:code-availability]{Code availability paragraph} for code that checks the computations in the following proofs.
\begin{theorem}\label{thm:Xa}
    For $a \in \bb Z$, let $X_a := Z''_{13}(a, 2, a+2, a+1)$. Then $X_a$ is a smooth non-projective proper toric threefold, and $\omega_{X_a}$ is $\otimes$-generating for every $a$.
\end{theorem}
\begin{proof}
    Non-projectivity holds by \Cref{rem:Z13-attribution}, since $b = 2 \neq 0$ and $c - a = 2 \neq 0$. Set $D := -K_{X_a} = \sum_{i=1}^8 D_{v_i}$ to be the anti-canonical divisor, which is Cartier and big since $X_a$ is smooth and by \Cref{cor:big-cone}. By \Cref{thm:implementable} and \Cref{rem:intermediate-cones}, it suffices to show that for each of the eight rays $\rho_i = \operatorname{cone}(v_i)$, one of $P_{\pm D}(\rho_i)$ is $2$-dimensional, and for each of the eighteen $2$-dimensional cones $\tau$, one of $P_{\pm D}(\tau)$ is $1$-dimensional.

    \emph{Rays.} By \Cref{cor:conv-aff}, it suffices to show that each $v_i$ is a vertex of the convex hull of $v_i$ and the generators sharing a cone with it. For $i=1,2$, the supporting plane $x=1$ meets this convex hull in $[v_1,v_2]$. For $i=5,8$, the supporting plane $x=-1$ meets it in $[v_5,v_8]$. For the remaining rays, the following supporting planes meet the corresponding convex hulls in the indicated edges:
    \[
    \begin{array}{c|c|c}
    \text{Ray} & \text{Supporting plane} & \text{Edge} \\\hline
    \rho_3 & x-2y+2z=2 & [v_3,v_7] \\
    \rho_4 & x-2y+2z=-2 & [v_4,v_6] \\
    \rho_6 & -x-2y+2z=-2 & [v_4,v_6] \\
    \rho_7 & -x-2y+2z=2 & [v_3,v_7]
    \end{array}
    \]
    Each $v_i$ is an endpoint of the corresponding edge, hence a vertex, for every $a \in \bb Z$.

    \emph{Walls.} Write $ij$ for the cone $\operatorname{cone}(v_i, v_j)$. Since every invariant curve is isomorphic to $\bb P^1$, it suffices to check that every wall has nonzero anticanonical degree. By \Cref{eq:wall-degree},
    \[
    D\cdot V(34)=2a+3,\qquad D\cdot V(67)=-2a-5,
    \]
    which are nonzero for every integer $a$. The remaining sixteen walls have nonzero anticanonical degree for every $a$:
    \[
    \begin{array}{c|r||c|r}
    \tau & D\cdot V(\tau) & \tau & D\cdot V(\tau) \\\hline
    12 & 2 & 37 & 2 \\
    14 & -1 & 38 & -1 \\
    16 & 2 & 45 & 2 \\
    17 & 3 & 46 & 4 \\
    23 & 2 & 56 & -1 \\
    24 & 3 & 58 & 2 \\
    27 & -1 & 68 & 3 \\
    35 & 3 & 78 & 2
    \end{array}
    \]
    Every wall and every ray therefore satisfies the criterion, and $\omega_{X_a} = \ecal O_{X_a}(-D)$ is $\otimes$-generating for all $a \in \bb Z$.
\end{proof}
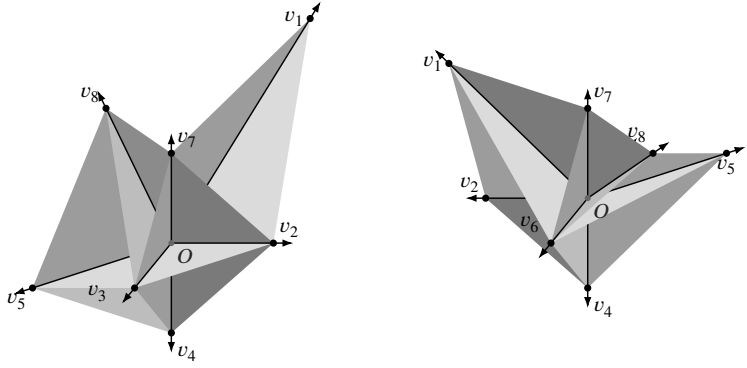
\begin{figure}[t] 
    \centering



%
%
%
%
%
\begin{tikzpicture}[
    x={(1.25cm,0.00cm)},
    y={(0.45cm,0.55cm)},
    z={(0.00cm,1.10cm)},
    cls plot, scale=1.08]
    \let\wallfill\CLSwall
    \let\fanray\CLSray
    \colorlet{g1}{clsgray1}
    \colorlet{g2}{clsgray2}
    \colorlet{g3}{clsgray3}
    \colorlet{g4}{clsgray4}

    \begin{scope}
        \coordinate (O)  at ( 0, 0, 0);
        \coordinate (v1) at ( 1, 1, 2);
        \coordinate (v2) at ( 1, 0, 0);
        \coordinate (v3) at ( 0,-1, 0);
        \coordinate (v4) at ( 0, 0,-1);
        \coordinate (v5) at (-1,-1, 0);
        \coordinate (v6) at ( 0, 1, 0);
        \coordinate (v7) at ( 0, 0, 1);
        \coordinate (v8) at (-1, 1, 1);
    
    
        \wallfill{v5}{v8}{g3} 
        \wallfill{v7}{v8}{g3!50!g4} 
        \wallfill{v1}{v7}{g3}   
        \wallfill{v1}{v2}{g1}  
        \wallfill{v2}{v4}{g4}   
        \wallfill{v4}{v5}{g2}   
    
        \fanray{v1}
        
        \wallfill{v3}{v5}{g1} 
        \wallfill{v3}{v4}{g3}   
    
        \fanray{v4}
        \fanray{v5}
    
        \wallfill{v3}{v8}{g2} 
    
        \fanray{v8}
        
        \wallfill{v2}{v7}{g4}   
        \wallfill{v3}{v7}{g3} 
        \wallfill{v2}{v3}{g1}   
    
        \fanray{v2}
        \fanray{v3}
        \fanray{v7}
        
        \node[origin] at (O) {};
        \foreach \v in {v1,v2,v3,v4,v5,v7,v8} {
            \node[gen] at (\v) {};
        }
    
        \node[lbl, anchor=east] at (v1) {$v_1$};
        \node[lbl, anchor=south west] at (v2) {$v_2$};
        \node[lbl, anchor=east, xshift=-8pt, yshift=-2pt] at (v3) {$v_3$};
        \node[lbl, anchor=north west, yshift=-3pt] at (v4) {$v_4$};
        \node[lbl, anchor=north east] at (v5) {$v_5$};
        \node[lbl, anchor=south west] at (v7) {$v_7$};
        \node[lbl, anchor=south east] at (v8) {$v_8$};
        \node[lbl, text=gray!30!black, anchor=north west] at (O) {$O$};
    \end{scope}

    \begin{scope}[xshift=5.1cm, yshift=5.5mm]
        \coordinate (O)  at ( 0, 0, 0);
        \coordinate (v1) at (-1,-1, 2);
        \coordinate (v2) at (-1, 0, 0);
        \coordinate (v3) at ( 0, 1, 0);
        \coordinate (v4) at ( 0, 0,-1);
        \coordinate (v5) at ( 1, 1, 0);
        \coordinate (v6) at ( 0,-1, 0);
        \coordinate (v7) at ( 0, 0, 1);
        \coordinate (v8) at ( 1,-1, 1);
    

        \fanray{v3}
    
        \wallfill{v1}{v7}{g4}
        \wallfill{v7}{v8}{g4} 
        \wallfill{v8}{v5}{g3}
        \wallfill{v5}{v4}{g3}
        \wallfill{v4}{v2}{g4}
        \wallfill{v2}{v1}{g3}

        \fanray{v2}
        
        \wallfill{v1}{v4}{g1}

        \wallfill{v6}{v1}{g1}
        \wallfill{v6}{v4}{g2}

        \fanray{v1}
        \fanray{v4}

        \wallfill{v6}{v7}{g3}
        \wallfill{v6}{v5}{g1}

        \fanray{v5}
        \fanray{v7}

        \wallfill{v6}{v8}{g2}

        \fanray{v6}
        \fanray{v8}
        
        \node[origin] at (O) {};
        \foreach \v in {v1,v6,v2,v4,v5,v7,v8} {
            \node[gen] at (\v) {};
        }
    
        \node[lbl, anchor=east] at (v1) {$v_1$};
        \node[lbl, anchor=south east] at (v2) {$v_2$};
        \node[lbl, anchor=north west, yshift=-3pt] at (v4) {$v_4$};
        \node[lbl, anchor=north] at (v5) {$v_5$};
        \node[lbl, anchor=south east, xshift=-2pt, yshift=2pt] at (v6) {$v_6$};
        \node[lbl, anchor=south west] at (v7) {$v_7$};
        \node[lbl, anchor=south east, yshift=2pt] at (v8) {$v_8$};
        \node[lbl, text=gray!30!black, anchor=north west] at (O) {$O$};
    
    \end{scope}
\end{tikzpicture}


    \caption{The fan of $X_{-1}=Z_{13}''(-1,2,1,0)$ as in \Cref{thm:Xa}.} \label{fig:X-1}
\end{figure}

\begin{theorem}\label{thm:Ya}
    For $a \in \bb Z$, let $Y_a := Z''_{13}(a, 2, a+2, a)$. Then $Y_a$ is a smooth non-projective proper toric threefold, and $\omega_{Y_a}$ is $\otimes$-generating if and only if $a \neq -1, -2$.
\end{theorem}
\begin{proof}
    Non-projectivity holds as before. Set $D := -K_{Y_a}$.

    \emph{Rays.} The planes in the proof of \Cref{thm:Xa} support the corresponding convex hulls for $Y_a$ and meet them in the same edges, so every ray satisfies the criterion for every $a \in \bb Z$.

    \emph{Walls.} By \Cref{eq:wall-degree},
    \[
    D\cdot V(34)=2a+2,\qquad D\cdot V(67)=-2a-4,
    \]
    so by \Cref{prop:wall-criterion}, the wall $34$ fails the criterion exactly at $a=-1$, and the wall $67$ exactly at $a=-2$. The remaining sixteen walls pass for every $a$: their anticanonical degrees are as in the proof of \Cref{thm:Xa}, except that $D\cdot V(37)=3$ and $D\cdot V(46)=3$.
    Hence $\omega_{Y_a}$ is $\otimes$-generating if and only if $a \neq -1, -2$, with the obstruction carried by the invariant curve $V(34)$ for $a = -1$ and $V(67)$ for $a = -2$.
\end{proof}

\begin{remark}\label{rem:Ya-obstruction}
    For $a = -1, -2$, the obstructing subvariety of $Y_a$ is the invariant curve $V(\tau)$, where $\tau$ is $\operatorname{cone}(v_3, v_4)$ for $a = -1$ and $\operatorname{cone}(v_6, v_7)$ for $a = -2$. In the first case, the adjacent rays are $v_2$ and $v_5$, and the wall relation $v_2 + v_5 = (0,a,a) = -a\,v_3 - a\,v_4$ has $(b_1, b_2) = (1,1)$ at $a = -1$. In the second case, the adjacent rays are $v_1$ and $v_8$, and $v_1 + v_8 = (0,a+3,a+3) = (a+3)v_6 + (a+3)v_7$ has $(b_1,b_2) = (1,1)$ at $a = -2$. In both cases $V(\tau)$ is a $(-1,-1)$-curve in the sense of \Cref{prop:wall-types-smooth} below and the corresponding Atiyah flop gives a (non-projective) Fourier--Mukai partner. In \cite{fujino2024non}*{Lemma 5.1}, Fujino and Sato give a sufficient condition on $(a,b,c,d)$ for the fan of $Z''_{13}$ to be the unique nonsingular complete fan on its rays, which implies that there is no Atiyah flop (or more generally a birational map that is an isomorphism in codimension $1$) to another (smooth complete toric) variety. On the other hand, their condition is not necessary as observed in \cite{fujino2024non}*{Remark 5.3}. Our subfamilies $X_a$ and $Y_a$ also lie outside the range of their sufficient condition for most $a$  and provide more instances of \cite{fujino2024non}*{Remark 5.3}.  More research on toric/combinatorial understanding of the birational rigidity coming from $\otimes$-generation of the canonical bundles and the $\otimes$-generating fan may be an interesting future direction.
\end{remark}

\begin{question}\label{q:Z13-general}
    For which $(a,b,c,d)$ is the canonical bundle of $Z''_{13}(a,b,c,d)$ $\otimes$-generating? In the subfamilies above, every obstruction is an invariant $(-1,-1)$-curve. Is that the case throughout the family?
\end{question}

\subsection{Blow-ups of \texorpdfstring{$\bb P^3$}{P3}}\label{subsec:P3-blowups}

Let $\bb P^3 = X_\Sigma$ for the fan $\Sigma$ in $\bb Z^3$ with rays generated by
\[
u_1 = e_1, \quad u_2 = e_2, \quad u_3 = e_3, \quad u_4 = -e_1 - e_2 - e_3.
\]
The six invariant lines $L_{ij}$ correspond to $\operatorname{cone}(u_i,u_j)$ for $1\leq i<j\leq4$, and blowing up $L_{ij}$ inserts $e_{ij}:=u_i+u_j$. For an \emph{ordered} subset $S$ of these lines, let $X_S$ be the iterated blow-up and let $P_S$ be the convex hull of $u_1,\dots,u_4$ and the $e_{ij}$ with $L_{ij}\in S$. \Cref{fig:P3-lines} illustrates the polytope and the proof below.

\begin{prop}\label{prop:P3-blowup-obstructions}
    The canonical bundle $\omega_{X_S}$ fails to be $\otimes$-generating if and only if there exists indices $i,j,k\in [4]$ for which $L_{ij}, L_{jk}\in S$ but $L_{ik}\notin S$. In this case, every obstructing invariant subvariety is a $(-1,-1)$-curve, inducing an Atiyah flop.
\end{prop}

\begin{proof}
  We apply the hyperplane criterion (\Cref{cor:conv-aff}) on $P_S$. First, one can check that every ray generator $u_i$ or $e_{ij}$ is a vertex of $P_S$, so it lies strictly outside of the convex hull of the other generators. Hence, no invariant surface can obstruct (see \Cref{fig:P3-lines} (b)). 

  For invariant curves, observe that any wall $\operatorname{cone}(w_1, w_2)$ whose generators $w_1, w_2$ span an edge of $P_S$ satisfies the hyperplane criterion. Indeed, a supporting hyperplane exposing the edge $\overline{w_1w_2}$ contains $\operatorname{aff}(w_1,w_2)$ and places all other vertices strictly on one side. Thus, for the opposite generators $x_1,x_2$, we have $\operatorname{aff}(w_1,w_2)\cap\operatorname{conv}(x_1,x_2)=\varnothing$. Furthermore, if a wall $\operatorname{cone}(w_1, w_2)$ is contained in the boundary of a maximal cone of the original $\mathbb P^3$, then $\operatorname{conv}(w_1, w_2)$ is an edge of $P_S$, so the above argument applies, which means such a curve cannot obstruct (see \Cref{fig:P3-lines} (c)).

  It remains to examine the interior walls. Fix an original cone $\operatorname{cone}(u_i, u_j, u_k)$ on which a blowup occurs. Relabel for convenience so that $L_{ij}$ is blown up first, introducing the wall $\operatorname{cone}(e_{ij}, u_k)$. Subsequent blowups at $L_{ik}, L_{jk}$, if present, introduce walls $\operatorname{cone}(e_{ij}, e_{ik})$ and $\operatorname{cone}(e_{ij}, e_{jk})$, but these are edges of $P_S$, so they do not obstruct. There are three cases we need to check.

  \emph{Case I.} If $L_{ij}\in S$ and $L_{ik}, L_{jk}\notin S$, then $\operatorname{conv}(e_{ij}, u_k)$ is also an edge of $P_S$, so there is no obstruction (see \Cref{fig:P3-lines} (d)).

  \emph{Case II.} This is the interesting case. If $L_{ij}, L_{ik}\in S$ and $L_{jk}\notin S$, then the two maximal cones are $\operatorname{cone}(e_{ij}, u_j, u_k)$ and $\operatorname{cone}(e_{ij}, e_{ik}, u_k)$ meeting along the wall $\operatorname{cone}(e_{ij}, u_k)$. We have the circuit \[e_{ij} + u_{k} = u_j+e_{ik},\] giving the obstruction \[\operatorname{conv}(e_{ik}, u_j) \cap \operatorname{aff}(e_{ij},u_k) = \frac{e_{ij} + u_{k}}{2} = \frac{u_j+e_{ik}}{2}.\]
  By \Cref{prop:wall-types-smooth}, the obstructing curve $V(\operatorname{cone}(e_{ij}, u_k))$ is of type $(-1,-1)$ (see \Cref{fig:P3-lines} (e)). 

  \emph{Case III.} If all of $L_{ij}, L_{ik}, L_{jk}\in S$, then the two maximal cones are $\operatorname{cone}(e_{ij}, e_{ik}, u_k)$ and $\operatorname{cone}(e_{ij}, e_{jk}, u_k)$ meeting along the wall $\operatorname{cone}(e_{ij}, u_k)$. These satisfy \[e_{ik} + e_{jk} = e_{ij} + 2u_k,\]
  which is no longer an obstruction (see \Cref{fig:P3-lines} (f)). This concludes the proof. 
\end{proof}

\begin{figure}[t]
    \centering
\begin{tikzpicture}[cls plot]
    \begin{scope}[shift={(-1.05,0)},scale=0.82]
        \coordinate (u1) at (-0.614466,1.671483);
        \coordinate (u2) at (-1.351825,-1.005357);
        \coordinate (u3) at (1.351825,0.104853);
        \coordinate (u4) at (0.614466,-0.770979);
        \coordinate (e12) at (-1.966291,0.666126);
        \coordinate (e13) at (0.737359,1.776336);
        \coordinate (e14) at (0.000000,0.900504);
        \coordinate (e23) at (0.000000,-0.900504);
        \coordinate (e24) at (-0.737359,-1.776336);
        \coordinate (e34) at (1.966291,-0.666126);
        \fill[clsface2] (e13) -- (e12) -- (e23) -- cycle;
        \fill[clsface1] (u1) -- (e12) -- (e13) -- cycle;
        \fill[clsface2] (e23) -- (e12) -- (u2) -- cycle;
        \fill[clsface2] (u3) -- (e13) -- (e23) -- cycle;
        \fill[clsface1] (e34) -- (e13) -- (u3) -- cycle;
        \fill[clsface1] (e34) -- (e23) -- (e24) -- cycle;
        \fill[clsface1] (e24) -- (e23) -- (u2) -- cycle;
        \fill[clsface1] (u3) -- (e23) -- (e34) -- cycle;
        \draw[edgeBase] (e12) -- (e13);
        \draw[edgeHidden] (e12) -- (e14);
        \draw[edgeBase] (e12) -- (e23);
        \draw[edgeHidden] (e12) -- (e24);
        \draw[edgeBase] (e12) -- (u1);
        \draw[edgeBase] (e12) -- (u2);
        \draw[edgeHidden] (e13) -- (e14);
        \draw[edgeBase] (e13) -- (e23);
        \draw[edgeBase] (e13) -- (e34);
        \draw[edgeBase] (e13) -- (u1);
        \draw[edgeBase] (e13) -- (u3);
        \draw[edgeHidden] (e14) -- (e24);
        \draw[edgeHidden] (e14) -- (e34);
        \draw[edgeHidden] (e14) -- (u1);
        \draw[edgeHidden] (e14) -- (u4);
        \draw[edgeBase] (e23) -- (e24);
        \draw[edgeBase] (e23) -- (e34);
        \draw[edgeBase] (e23) -- (u2);
        \draw[edgeBase] (e23) -- (u3);
        \draw[edgeBase] (e24) -- (e34);
        \draw[edgeBase] (e24) -- (u2);
        \draw[edgeHidden] (e24) -- (u4);
        \draw[edgeBase] (e34) -- (u3);
        \draw[edgeHidden] (e34) -- (u4);
        \foreach \v in {u1,u2,u3,u4,e12,e13,e14,e23,e24,e34} \node[gen] at (\v) {};
        \node[lbl,anchor=south] at (u1) {$u_{1}$};
        \node[lbl,anchor=north east] at (u2) {$u_{2}$};
        \node[lbl,anchor=south west] at (u3) {$u_{3}$};
        \node[lbl,anchor=north] at (u4) {$u_{4}$};
        \node[lbl,anchor=east] at (e12) {$e_{12}$};
        \node[lbl,anchor=west] at (e13) {$e_{13}$};
        \node[lbl,anchor=south] at (e14) {$e_{14}$};
        \node[lbl,anchor=north] at (e23) {$e_{23}$};
        \node[lbl,anchor=east] at (e24) {$e_{24}$};
        \node[lbl,anchor=west] at (e34) {$e_{34}$};
        \node[lbl,align=center,anchor=north] at (0,-2.65) {(a) The polytope $P_S$ where $S$ contains all 6 lines.};
    \end{scope}
    \begin{scope}[shift={(3.575,0)},scale=0.82]
        \coordinate (u1) at (-0.614466,1.690415);
        \coordinate (u2) at (-1.351825,-1.186441);
        \coordinate (u3) at (1.351825,-0.346483);
        \coordinate (u4) at (0.614466,-0.157492);
        \coordinate (e12) at (-1.966291,0.503975);
        \coordinate (e13) at (0.737359,1.343933);
        \coordinate (e14) at (0.000000,1.532923);
        \coordinate (e23) at (0.000000,-1.532923);
        \coordinate (e24) at (-0.737359,-1.343933);
        \coordinate (e34) at (1.966291,-0.503975);
        \fill[clsgray2] (e14) -- (e12) -- (e13) -- cycle;
        \fill[clsface2] (e13) -- (e12) -- (e23) -- cycle;
        \fill[clsface2] (e23) -- (e12) -- (u2) -- cycle;
        \fill[clsface2] (u3) -- (e13) -- (e23) -- cycle;
        \fill[clsface1] (e34) -- (e13) -- (u3) -- cycle;
        \fill[clsface1] (u3) -- (e23) -- (e34) -- cycle;
        \draw[edgeBase] (e12) -- (e13);
        \draw[edgeBase] (e12) -- (e14);
        \draw[edgeBase] (e12) -- (e23);
        \draw[edgeHidden] (e12) -- (e24);
        \draw[edgeBase] (e12) -- (u2);
        \draw[edgeBase] (e13) -- (e14);
        \draw[edgeBase] (e13) -- (e23);
        \draw[edgeBase] (e13) -- (e34);
        \draw[edgeBase] (e13) -- (u3);
        \draw[edgeHidden] (e14) -- (e24);
        \draw[edgeHidden] (e14) -- (e34);
        \draw[edgeHidden] (e14) -- (u4);
        \draw[edgeHidden] (e23) -- (e24);
        \draw[edgeBase] (e23) -- (e34);
        \draw[edgeBase] (e23) -- (u2);
        \draw[edgeBase] (e23) -- (u3);
        \draw[edgeHidden] (e24) -- (e34);
        \draw[edgeHidden] (e24) -- (u2);
        \draw[edgeHidden] (e24) -- (u4);
        \draw[edgeBase] (e34) -- (u3);
        \draw[edgeHidden] (e34) -- (u4);
        \draw[edgeHidden] (u1) -- (e12);
        \draw[edgeHidden] (u1) -- (e13);
        \draw[edgeHidden] (u1) -- (e14);
        \foreach \v in {u2,u3,u4,e12,e13,e14,e23,e24,e34} \node[gen] at (\v) {};
        \node[gen,inner sep=1.6pt] at (u1) {};
        \node[lbl,anchor=south] at (u1) {$u_i$};
        \node[lbl,anchor=east] at (e12) {$e_{ij}$};
        \node[lbl,anchor=west] at (e13) {$e_{ik}$};
        \node[lbl,anchor=north] at (e14) {$e_{il}$};
        \node[lbl,align=center,anchor=north] at (0,-2.65) {(b) No surface obstructs.};
    \end{scope}
    \begin{scope}[shift={(8.2,0)},scale=0.82]
        \coordinate (u1) at (-0.285774,0.885949);
        \coordinate (u2) at (-0.842666,0.396051);
        \coordinate (u3) at (1.985761,0.396051);
        \coordinate (u4) at (0.285774,-1.750051);
        \coordinate (e12) at (-1.414214,1.300000);
        \coordinate (e13) at (1.414214,1.300000);
        \coordinate (plane0) at (-1.885000,1.120000);
        \coordinate (plane1) at (2.515000,1.120000);
        \coordinate (plane2) at (1.885000,1.480000);
        \coordinate (plane3) at (-2.515000,1.480000);
        \fill[clsface2] (plane0) -- (plane1) -- (plane2) -- (plane3) -- cycle;
        \draw[edgeBase] (plane0) -- (plane1) -- (plane2) -- (plane3) -- cycle;
        \fill[clsgray2] (e13) -- (e12) -- (u2) -- (u3) -- cycle;
        \fill[clsface1] (u2) -- (e12) -- (u4) -- cycle;
        \fill[clsface1] (u3) -- (u2) -- (u4) -- cycle;
        \draw[edgeBase] (e12) -- (e13);
        \draw[edgeHidden] (e12) -- (u1);
        \draw[edgeBase] (e12) -- (u2);
        \draw[edgeBase] (e12) -- (u4);
        \draw[edgeHidden] (e13) -- (u1);
        \draw[edgeBase] (e13) -- (u3);
        \draw[edgeHidden] (e13) -- (u4);
        \draw[edgeHidden] (u1) -- (u4);
        \draw[edgeBase] (u2) -- (u3);
        \draw[edgeBase] (u2) -- (u4);
        \draw[edgeBase] (u3) -- (u4);
        \draw[edgeBlow] ($(e12)!-0.22!(e13)$) -- ($(e13)!-0.22!(e12)$);
        \draw[chord,line width=0.85pt] (u1) -- (u3);
        \foreach \v in {u1,u2,u3,u4,e12,e13} \node[gen] at (\v) {};
        \node[lbl,anchor=south east] at (e12) {$w_1$};
        \node[lbl,anchor=south west] at (e13) {$w_2$};
        \node[lbl,anchor=north east] at (u1) {$x_1$};
        \node[lbl,anchor=north west] at (u3) {$x_2$};
        \node[lbl] at (0,1.85) {supporting plane};
        \node[lbl] at (0.35,-1.05) {$P_S$};
        \node[lbl,align=center,anchor=north] at (0,-2.65) {(c) An edge of $P_S$ does not obstruct.};
    \end{scope}
    \begin{scope}[shift={(-1.05,-5.15)},scale=0.82]
        \coordinate (ui) at (-1.811250,-1.748000);
        \coordinate (uj) at (0.948750,-0.920000);
        \coordinate (uk) at (1.811250,0.920000);
        \coordinate (eij) at (-1.811250,-0.920000);
        \fill[clsface2] (uj) -- (eij) -- (ui) -- cycle;
        \fill[clsface1] (uk) -- (eij) -- (uj) -- cycle;
        \draw[edgeBase] (ui) -- (uj);
        \draw[edgeHidden] (ui) -- (uk);
        \draw[edgeBase] (ui) -- (eij);
        \draw[edgeBase] (uj) -- (uk);
        \draw[edgeBase] (uj) -- (eij);
        \draw[edgeBase] (uk) -- (eij);
        \draw[edgeBlow] ($(eij)!-0.13!(uk)$) -- ($(uk)!-0.13!(eij)$);
        \draw[chord,line width=0.85pt] (ui) -- (uj);
        \foreach \v in {ui,uj,uk,eij} \node[gen] at (\v) {};
        \node[lbl,anchor=north east] at (ui) {$u_i$};
        \node[lbl,anchor=north west] at (uj) {$u_j$};
        \node[lbl,anchor=south east] at (eij) {$e_{ij}$};
        \node[lbl,anchor=south east] at (uk) {$u_k$};
        \node[lbl,align=center,anchor=north] at (0,-2.65) {(d) Case I: $\{L_{ij}\}$.};
    \end{scope}
    \begin{scope}[shift={(3.575,-5.15)},scale=0.82]
        \coordinate (eij) at (-1.811250,-0.920000);
        \coordinate (uj) at (0.948750,-0.920000);
        \coordinate (uk) at (1.811250,0.920000);
        \coordinate (eik) at (-0.948750,0.920000);
        \fill[clsgray2] (eij) -- (uj) -- (uk) -- (eik) -- cycle;
        \draw[edgeBase] (eij) -- (uj) -- (uk) -- (eik) -- cycle;
        \draw[edgeBlow] ($(eij)!-0.13!(uk)$) -- ($(uk)!-0.13!(eij)$);
        \draw[chord,line width=0.85pt] (uj) -- (eik);
        \foreach \v in {eij,uj,uk,eik} \node[gen] at (\v) {};
        \node[gen,inner sep=1.5pt] at ($(eij)!0.5!(uk)$) {};
        \node[lbl,anchor=south east] at (eij) {$e_{ij}$};
        \node[lbl,anchor=north west] at (uj) {$u_j$};
        \node[lbl,anchor=south east] at (eik) {$e_{ik}$};
        \node[lbl,anchor=south east] at (uk) {$u_k$};
        \node[lbl,align=center,anchor=north] at (0,-2.65) {(e) Case II: $\{L_{ij},L_{ik}\}$.};
    \end{scope}
    \begin{scope}[shift={(8.2,-5.15)},scale=0.82]
        \coordinate (eij) at (-1.811250,-0.920000);
        \coordinate (eik) at (-0.948750,0.920000);
        \coordinate (ejk) at (1.811250,1.748000);
        \coordinate (uk) at (1.811250,0.920000);
        \fill[clsgray2] (eik) -- (eij) -- (ejk) -- cycle;
        \draw[edgeBase] (eij) -- (eik);
        \draw[edgeBase] (eij) -- (ejk);
        \draw[edgeBase] (eij) -- (uk);
        \draw[edgeBase] (eik) -- (ejk);
        \draw[edgeHidden] (eik) -- (uk);
        \draw[edgeBase] (ejk) -- (uk);
        \draw[edgeBlow] ($(eij)!-0.13!(uk)$) -- ($(uk)!-0.13!(eij)$);
        \draw[chord,line width=0.85pt] (eik) -- (ejk);
        \foreach \v in {eij,eik,ejk,uk} \node[gen] at (\v) {};
        \node[lbl,anchor=south east] at (eij) {$e_{ij}$};
        \node[lbl,anchor=south east] at (eik) {$e_{ik}$};
        \node[lbl,anchor=south west] at (ejk) {$e_{jk}$};
        \node[lbl,anchor=north west] at (uk) {$u_k$};
        \node[lbl,align=center,anchor=north] at (0,-2.65) {(f) Case III: $\{L_{ij},L_{ik},L_{jk}\}$.};
    \end{scope}
    \draw[edgeBlow] (-0.5,-8.25) -- (0.25,-8.25);
    \node[lbl,anchor=west] at (0.4,-8.25) {affine line through the wall generators};
    \draw[chord,line width=0.85pt] (-0.5,-8.7) -- (0.25,-8.7);
    \node[lbl,anchor=west] at (0.4,-8.7) {segment joining the opposite generators};
\end{tikzpicture}

    \caption{An illustration of the polytope $P_S$ and the separation arguments used in \Cref{prop:P3-blowup-obstructions}. In (d)--(f), the subsets shown are $S\cap\{L_{ij},L_{ik},L_{jk}\}$, with $L_{ij}$ blown up first. Only Case II obstructs.}
    \label{fig:P3-lines}
\end{figure}

\begin{remark}\label{rem:non-monotone}
    The above classification exhibits non-monotonicity of $\otimes$-generation of canonical bundles under further blowups. Every $2$-element subset of a concurrent triple is an intersecting pair, so blowing up two of three concurrent lines obstructs $\otimes$-generation while blowing up the third removes the obstruction. Likewise, the blowup along the union of all six lines is $\otimes$-generating although many subsets in between are not.
\end{remark}

\begin{remark}
    One can also check that the blowup of $\mathbb P^3$ along invariant lines in $S$ has $\otimes$-generating canonical bundle if and only if every facet of the polytope $P_S$ is a triangle. Conversely, this blowup does not have $\otimes$-generating canonical bundle if and only if the polytope $P_S$ has a face that is a parallelogram. 
\end{remark}

\subsection{Projective bundles}
In contrast to products of complete Gorenstein toric varieties, where the canonical bundle is $\otimes$-generating if and only if it is anti-ample, toric projective bundles can admit canonical bundles that are $\otimes$-generating, but not anti-ample (nor ample). Although we should be able to give a complicated formula for when projective bundles have $\otimes$-generating canonical bundles, we focus on the following simple case for clarity and we use the derived Kleiman criterion (\Cref{thm:curve-criterion}). Note that the proof is quite analogous to ones for examples in \cite{itoolander2026botamestack}*{Section 2.3}. 
\begin{prop}\label{prop: projective bundle}
    Let $X$ be a proper Gorenstein toric variety with $\dim X >0$ and let $\ecal L$ be a line bundle. Set $\pi: Y = \bb P_X(\ecal O_X^{\oplus n} \oplus \ecal L) \to X$ for $n \geq 1$. If $\omega_X \otimes \ecal L$ is anti-ample and $\omega_X \otimes \ecal L^{\otimes - n}$ is $\otimes$-generating, then $\omega_Y$ is $\otimes$-generating. The converse holds if $n \geq 2$. 
\end{prop}
\begin{proof}
    First of all, we have
    \[
    \omega_Y \cong \ecal O_Y(-n-1) \otimes \pi^*(\omega_X \otimes \det (\ecal O_X^{\oplus n} \oplus \ecal L)) \cong  \ecal O_Y(-n-1) \otimes \pi^*(\omega_X \otimes  \ecal L).
    \]
    First, suppose that $\omega_X \otimes \ecal L$ is anti-ample and $\omega_X \otimes \ecal L^{\otimes - n}$ is $\otimes$-generating. Take an irreducible complete curve $C \subset Y$ and it suffices to show $\deg_C \omega_Y|_C \neq 0$ by \Cref{thm:curve-criterion}. Set $f: C \subset Y \to X$ to be the composition and write $\ecal Q = \ecal O_Y(1)|_{C}$ so that $\omega_Y|_C \cong \ecal Q^{\otimes(-n-1)} \otimes f^*(\omega_X \otimes \ecal L)$. Also, the morphism $f$ corresponds to the quotient 
    \[
    f^*\ecal L \oplus \ecal O_C^{\oplus n} \to \ecal Q.
    \]
    First, if $f$ is constant, then $C$ is contained in a fiber of $\pi$ and hence $\omega_Y|_C \cong \ecal Q^{\otimes(-n-1)}$ is anti-ample and we are done. Now, suppose $f$ is finite. If there is a component so that $\ecal O_C \to f^*\ecal L \oplus \ecal O_C^{\oplus n} \to \ecal Q$ is non-zero, then $\deg \ecal Q \geq 0$, and since $f^*(\omega_X \otimes \ecal L)$ is anti-ample, $\deg_C \omega_Y|_C < 0$. Finally, if all the components $\ecal O_C \to f^*\ecal L \oplus \ecal O_C^{\oplus n} \to \ecal Q$ are zero, then $f^* \ecal L$ surjects onto $\ecal Q$ and hence $f^*\ecal L \cong \ecal Q$. In particular, 
    \[
    \omega_Y|_C \cong f^*(\ecal L^{\otimes (-n-1)} \otimes \omega_X \otimes \ecal L) = f^*(\ecal L^{\otimes - n} \otimes \omega_X)
    \]
    is $\otimes$-generating by assumption together with \cite{ito2025derived}*{Lemma 2.4}, as $f$ is finite. Thus, $\deg \omega_Y|_C \neq  0$ by \cite{ito2025derived}*{Proposition 5.1}. 
    
    Conversely, suppose that $\omega_Y$ is $\otimes$-generating and that $n \geq 2$. First, consider the section $s: X \to Y$ corresponding to the quotient $\ecal O_X^{\oplus n} \oplus \ecal L \to \ecal L$. Then, 
    \[
    s^* \omega_Y \cong \ecal L^{\otimes(-n-1)} \otimes \omega_X \otimes \ecal L = \ecal L^{\otimes - n} \otimes \omega_X
    \]
    is $\otimes$-generating. Similarly, consider the map $f: \bb P_X (\ecal O_X^{\oplus n}) \to Y$ of projective bundles corresponding to the quotient $\ecal O_X^{\oplus n} \oplus \ecal L \to \ecal O_X^{\oplus n}$. Then
    \[
    f^* \omega_Y \cong \ecal O_{\bb P_X(\ecal O_X^{\oplus n})}(-n-1) \otimes (\bb P_X (\ecal O_X^{\oplus n}) \to X)^* (\omega_X \otimes \ecal L). 
    \]
    Therefore, under the identification $X \times \bb P^{n-1} \cong \bb P_X (\ecal O_X^{\oplus n})$, $f^* \omega_Y$ pulls back to $(\omega_X \otimes \ecal L) \boxtimes \ecal O_{\bb P^{n-1}}(-n-1)$ and since it is $\otimes$-generating and the latter component is anti-ample, $\omega_X \otimes \ecal L$ is anti-ample by \Cref{thm:products}. 
\end{proof}
Let us apply this criterion over the Hirzebruch surface $\bb F_m$ for $m \geq 0$.
\begin{example}
Let $S \subset \bb F_m = X$ be the section with $S^2 = -m$ and let $F$ be a fiber divisor. Recall that any line bundle $\ecal L$ on $X$ can be written as $\ecal L = \ecal O_{X}(aS +bF)$. It is easy to see $\ecal L$ is ample if and only if $a>0$ and $b>am$. Hence, 
\[
\omega_{X} \otimes \ecal L = \ecal O_X((a-2) S + (b-m-2) F)  \text{ and } \omega_X \otimes \ecal L^{\otimes -n} = \ecal O_X((-na -2)S + (-nb-m-2) F ).
\]
Set $Y = \bb P_{X}(\ecal O_{X}^{\oplus n} \oplus \ecal L)$. First, $\omega_X\otimes \ecal L$ is anti-ample if and only if $$a - 2< 0 \quad \text{and} \quad b-m-2 < (a-2)m.$$ Recall by \cite{ito2025derived}*{Lemma 5.3} that $\omega_X \otimes \ecal L^{\otimes -n}$ is $\otimes$-generating if and only if $$(-na-2)(-nb-m-2) > 0$$ (so that it is big or anti-big) and $$((\omega_X \otimes \ecal L^{\otimes -n}) . S) = -m(-na-2) + (-nb -m -2) \neq 0.$$
It is straightforward to observe the last condition is redundant if $n \geq 2$, so in this case, we see that $\omega_{X} \otimes \ecal L$ is anti-ample and $\omega_X \otimes \ecal L^{\otimes -n}$ is $\otimes$-generating if and only if 
\[
a<2, \quad b-am < 2-m, \quad (na+2)(nb+m +2)>0. 
\]
For example, the canonical bundle $\omega_Y$ of
\[
Y = \bb P_X(\ecal O_X^{\oplus n} \oplus \ecal O(S))
\]
is $\otimes$-generating for any $n \geq 2$ and $m \geq 0$, and if moreover $m \geq 1$, then $\omega_Y$ is neither nef nor anti-nef, since $\omega_Y$ pulls back to the anti-ample line bundle $\omega_{X} \otimes \ecal L$, while $\omega_Y$ also pulls back to $\omega_X \otimes \ecal L^{\otimes-n} \cong \ecal O_X(-(n+2)S - (m+2) F)$, which intersects positively with $S$. 
\end{example}

\begin{figure}[t]
    \centering
\begin{tikzpicture}[
    cls plot, scale=1.12,
    x={(1.25cm,0.05cm)},
    y={(-0.40cm,0.85cm)},
    z={(0,1.30cm)}]
    \let\wall\CLSwall
    \let\ray\CLSray
    \colorlet{gray1}{clsgray4}
    \colorlet{gray2}{clsgray3}
    \colorlet{gray3}{clsgray2}
    \colorlet{gray4}{clsgray1}

    \coordinate (O)  at ( 0, 0, 0);
    \coordinate (xp) at ( 1, 0, 0);
    \coordinate (xm) at (-1, 0, 1);
    \coordinate (yp) at (-2, 1, 0);
    \coordinate (ym) at (-2,-1, 0);
    \coordinate (zp) at ( 0, 0, 1);
    \coordinate (zm) at ( 0, 0,-1);

    \wall{xm}{yp}{gray1}
    \wall{zm}{yp}{gray2}
    \ray{yp}
    \wall{xm}{zp}{gray2}
    \wall{xm}{ym}{gray3}
    \ray{xm}
    \wall{zm}{ym}{gray1}
    \wall{zm}{xp}{gray2}
    \ray{zm}
    \wall{zp}{ym}{gray4}
    \wall{zp}{xp}{gray2}
    \wall{xp}{ym}{gray3}
    \ray{xp}\ray{ym}\ray{zp}

    \node[origin] at (O) {};
    \foreach \v in {xp,xm,yp,ym,zp,zm} {\node[gen] at (\v) {};}

    \node[lbl, anchor=north west] at (xp) {$(1,0,0)$};
    \node[lbl, anchor=south east, yshift=3pt] at (xm) {$(-1,0,1)$};
    \node[lbl, anchor=north east] at (yp) {$(-2,1,0)$};
    \node[lbl, anchor=north east, yshift=-3pt] at (ym) {$(-2,-1,0)$};
    \node[lbl, anchor=south west] at (zp) {$(0,0,1)$};
    \node[lbl, anchor=north west] at (zm) {$(0,0,-1)$};
    \node[lbl, text=gray!30!black, anchor=north west] at (O) {$O$};
\end{tikzpicture}

    \caption{The fan of the projective bundle $\bb P(\ecal{O}\oplus \ecal{O}(-C))$ over $F_4$, where $C$ is the $-4$ curve.}\label{fig:bundle-over-F4}
\end{figure}
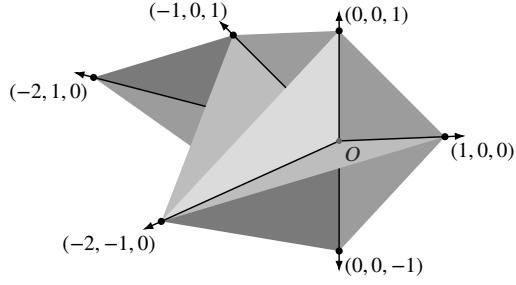

We can also directly observe this using \Cref{prop:hyperplane-test}.
\begin{example}\label{ex:bundle-over-F4}
    We can also see that on a projective bundle, $\otimes$-generation plus bigness won't imply ampleness. Consider the Hirzebruch surface $\bb F_4$ associated with the complete fan $\Sigma_4'$ in $\bb R^2$ with rays
    \[
    \rho_1 = \operatorname{cone}(-1,0), \quad \rho_2 = \operatorname{cone}(-2,1), \quad \rho_3 = \operatorname{cone}(1,0), \quad \rho_4 = \operatorname{cone}(-2,-1),
    \]
    so that $S:= D_{\rho_1}$ is the section of self-intersection $-4$. Let $Y =\bb P_{\bb F_4}(\ecal O\oplus \ecal{O}(-S))$ be the projective completion of the line bundle $\ecal{O}(-S)$, whose fan is depicted in \Cref{fig:bundle-over-F4}. Using \Cref{prop:hyperplane-test}, we can see that $\omega_Y$ is $\otimes$-generating, and in fact its restriction to any fiber is $\otimes$-generating. However, $\omega_Y$ is neither ample nor anti-ample. 
\end{example}

\section{Obstruction for \texorpdfstring{$\otimes$}{tensor}-generation and nonstandard autoequivalences}\label{sec:obstruction}
In the previous sections, we have focused on $\otimes$-generating line bundles. In this section, we focus on understanding precisely how $\otimes$-generation of canonical bundles can fail. For convenience, let us introduce the following terminology. 
\begin{definition}
Let $X$ be a proper Gorenstein variety. We say that a subvariety $Z \subset X$ is \textbf{obstructing ($\otimes$-generation)} if $\omega_X|_Z$ is neither big nor anti-big. In particular, if $X$ is toric, $K_X$ is $\otimes$-generating if and only if there is no obstructing $T$-invariant subvariety.    
\end{definition}
By \cite{ito2025derived}*{Lemma 5.7}, we have a complete classification of obstructing subvarieties for smooth projective toric surfaces.
\begin{example}
    Let $X$ be a smooth projective toric surface over an algebraically closed field. Then, a subvariety is obstructing if and only if it is a $(-2)$-curve. 
\end{example}     
In the rest of the section, we first recognize obstructing invariant curves in higher dimension and then consider implications of obstructing invariant subvarieties to autoequivalences of derived categories. 
\subsection{Local combinatorial structure of obstructing invariant curves}\label{subsec:wall-classification}

In this subsection, we classify the local structure of obstructing invariant curves for simplicial toric varieties, building on the local obstructions we have seen in \Cref{subsec:P3-blowups}. Let $X_\Sigma$ be an $n$-dimensional simplicial proper toric variety and $\tau = \operatorname{cone}(u_1, \dots, u_{n-1}) \in \Sigma$ a cone of dimension $n-1$. Since $\Sigma$ is complete and simplicial, $\tau$ is a facet of exactly two maximal cones $\sigma_1 = \tau + \operatorname{cone}(v_1)$ and $\sigma_2 = \tau + \operatorname{cone}(v_2)$. The $n+1$ vectors $u_1, \dots, u_{n-1}, v_1, v_2$ span $N_{\bb R}$ and satisfy, up to scaling, a unique linear relation
\begin{equation}\label{eq:wall-relation}
    c_1 v_1 + c_2 v_2 = \sum_{i=1}^{n-1} b_i u_i, \qquad c_1, c_2 > 0,
\end{equation}
where positivity of $c_1, c_2$ holds because $v_1$ and $v_2$ lie on opposite sides of the hyperplane spanned by $\tau$.

\begin{prop}\label{prop:wall-criterion}
    In the setting above, suppose $X_\Sigma$ is $\bb Q$-Gorenstein. The invariant curve $V(\tau)$ is obstructing if and only if
    \[
    \operatorname{conv}(v_1, v_2) \cap \operatorname{aff}(u_1, \dots, u_{n-1}) \neq \emp,
    \]
    if and only if the relation \Cref{eq:wall-relation}, which is unique up to scaling, satisfies $c_1 + c_2 = \sum_i b_i$.
\end{prop}
\begin{proof}
    The star fan of $\tau$ is a complete one-dimensional fan, so $V(\tau) \iso \bb P^1$, and $V(\tau)$ is obstructing if and only if $\deg(\ecal O_X(lK_X)|_{V(\tau)}) = 0$ for $l > 0$ with $lK_X$ Cartier, if and only if the polytope condition of \Cref{thm:implementable} fails at $\tau$ for $D = -K_X$, by \Cref{lem:toric-bigness,lem:restriction-polytope}. The divisor $-K_X$ has all coefficients $1$, so its marking is the marking by primitive generators. By the argument of \Cref{prop:hyperplane-test} at the single cone $\tau$, the polytope condition at $\tau$ holds if and only if some affine hyperplane through $u_1, \dots, u_{n-1}$ has $v_1$ and $v_2$ strictly on one side. By the separation argument of \Cref{cor:conv-aff}, such a hyperplane exists if and only if $\operatorname{conv}(v_1, v_2)$ is disjoint from $\operatorname{aff}(u_1, \dots, u_{n-1})$. This proves the first equivalence.

    For the second, suppose $p \in \operatorname{conv}(v_1, v_2) \cap \operatorname{aff}(u_1, \dots, u_{n-1})$. Since $\sigma_1$ and $\sigma_2$ are simplicial, neither $v_1$ nor $v_2$ lies in $\operatorname{span}(\tau) \supseteq \operatorname{aff}(u_1, \dots, u_{n-1})$, so $p = \frac{c_1 v_1 + c_2 v_2}{c_1 + c_2}$ with $c_1, c_2 > 0$. Writing $p = \sum_i t_i u_i$ with $\sum_i t_i = 1$ produces the relation $c_1 v_1 + c_2 v_2 = (c_1 + c_2)\sum_i t_i u_i$. The $n+1$ vectors $u_1, \dots, u_{n-1}, v_1, v_2$ span $N_{\bb R}$, so their relation space is one-dimensional, and this relation is a scaling of \Cref{eq:wall-relation}. Its coefficients satisfy $\sum_i (c_1+c_2)t_i = c_1 + c_2$, and the identity $c_1 + c_2 = \sum_i b_i$ is invariant under scaling, so it holds for \Cref{eq:wall-relation}. Conversely, if $c_1 + c_2 = \sum_i b_i$, then $p := \frac{c_1v_1 + c_2v_2}{c_1+c_2} = \frac{\sum_i b_i u_i}{\sum_i b_i}$ lies in the intersection.
\end{proof}

When $V(\tau)$ is obstructing, we normalize \Cref{eq:wall-relation} by $c_1 + c_2 = \sum_i b_i$ and set
\[
p := \frac{c_1 v_1 + c_2 v_2}{c_1 + c_2} = \frac{\sum_i b_i u_i}{\sum_i b_i} \ \in \ \operatorname{conv}(v_1,v_2) \cap \operatorname{aff}(u_1, \dots, u_{n-1}),
\]
and let $T := \operatorname{conv}(u_1, \dots, u_{n-1})$ be the simplex spanned by the generators of $\tau$. The position of $p$ relative to $T$ stratifies the obstructing walls of \Cref{prop:wall-criterion}:
\begin{enumerate}
    \item $p \in \operatorname{relint}(T)$, equivalently $b_i > 0$ for all $i$;
    \item $p$ lies on a proper face of $T$, equivalently $b_i \geq 0$ for all $i$ with some $b_j = 0$;
    \item $p \in \operatorname{aff}(T) \setminus T$, equivalently $b_i < 0$ for some $i$.
\end{enumerate}

In the smooth threefold case, the trichotomy is the classification of $K$-trivial walls by the normal bundle of the wall curve.

\begin{prop}\label{prop:wall-types-smooth}
    Let $X = X_\Sigma$ be a smooth proper toric threefold and $\tau = \operatorname{cone}(u_1, u_2)$ a wall with adjacent rays $v_1, v_2$ as above. Then $c_1 = c_2 = 1$ in \Cref{eq:wall-relation}, so the wall relation becomes
    \begin{equation}\label{eq:smooth-wall-relation}
    v_1 + v_2 = b_1 u_1 + b_2 u_2, \qquad b_1, b_2 \in \bb Z,
    \end{equation}
    and the curve $C := V(\tau) \iso \bb P^1$ satisfies
    \[
    C \cdot D_{u_i} = -b_i, \qquad
    \ecal N_{C/X} \iso \ecal O_{\bb P^1}(-b_1) \oplus \ecal O_{\bb P^1}(-b_2), \qquad
    \deg(\omega_X|_C) = b_1 + b_2 - 2.
    \]
    In particular, $C$ is obstructing if and only if $b_1 + b_2 = 2$, if and only if
    \[
    \ecal N_{C/X} \iso \ecal O(-1)\oplus\ecal O(-1), \quad \ecal O \oplus \ecal O(-2), \quad \text{or} \quad \ecal O(k) \oplus \ecal O(-2-k) \ \text{with } k \geq 1,
    \]
    corresponding to the cases (i), (ii), (iii) of the trichotomy.
\end{prop}
\begin{proof}
    Since $\sigma_1$ is a smooth cone, $(u_1, u_2, v_1)$ is a basis of $N$, and we may write $v_2 = \alpha v_1 + \beta_1 u_1 + \beta_2 u_2$ with integer coefficients. Since $(u_1, u_2, v_2)$ is also a basis, the change-of-basis matrix, which is triangular with diagonal $(1, 1, \alpha)$, has determinant $\alpha = \pm 1$. Since $v_1$ and $v_2$ lie on opposite sides of the plane $\operatorname{span}(\tau)$, we have $\alpha = -1$, which gives \Cref{eq:smooth-wall-relation} with $b_i = \beta_i$.

    For the degrees, we use \Cref{lem:restriction-polytope}. The star fan $\operatorname{Star}(\tau)$ is the complete fan in $N(\tau)_{\bb R} \iso \bb R$ whose rays are generated by the images $\bar v_1, \bar v_2$ of $v_1, v_2$, and \Cref{eq:smooth-wall-relation} gives $\bar v_2 = -\bar v_1$, so $C = V(\tau) \iso \bb P^1$. Let $D = \sum_\rho a_\rho D_\rho$ be any invariant Cartier divisor, and choose Cartier data with $m_\tau = m_{\sigma_1}$. By the proof of \Cref{lem:restriction-polytope}, $\ecal O_X(D)|_C \iso \ecal O_C(D_\tau)$, where $D_\tau$ has coefficient $-\bra{m_{\sigma_1} - m_\tau, \bar v_1}' = 0$ along $\bar v_1$ and coefficient $-\bra{m_{\sigma_2} - m_{\sigma_1}, v_2}$ along $\bar v_2$, so
    \[
    \deg\bigl(\ecal O_X(D)|_C\bigr)
    = \bra{m_{\sigma_1}, v_2} - \bra{m_{\sigma_2}, v_2}
    = \bra{m_{\sigma_1}, v_2} + a_{v_2},
    \]
    where we used $\bra{m_{\sigma_2}, v_2} = -a_{v_2}$. Expanding $v_2 = -v_1 + b_1 u_1 + b_2 u_2$ against $\bra{m_{\sigma_1}, v_1} = -a_{v_1}$ and $\bra{m_{\sigma_1}, u_i} = -a_{u_i}$ yields
    \begin{equation}\label{eq:wall-degree}
    \deg\bigl(\ecal O_X(D)|_C\bigr) = a_{v_1} + a_{v_2} - b_1 a_{u_1} - b_2 a_{u_2}.
    \end{equation}
    Taking $D = D_{u_i}$ in \Cref{eq:wall-degree} gives $C \cdot D_{u_i} = -b_i$, and taking $D = -K_X = \sum_\rho D_\rho$, with all $a_\rho = 1$, gives $\deg(\omega_X^{-1}|_C) = 2 - b_1 - b_2$, that is, $\deg(\omega_X|_C) = b_1 + b_2 - 2$. Since $C = D_{u_1} \cap D_{u_2}$ is a transverse intersection of two smooth invariant divisors, the normal bundle is
    \[
    \ecal N_{C/X} \iso \ecal O_X(D_{u_1})|_C \oplus \ecal O_X(D_{u_2})|_C \iso \ecal O(-b_1) \oplus \ecal O(-b_2).
    \]
    The curve $C$ is obstructing if and only if $\deg(\omega_X|_C) = 0$, since a line bundle on $\bb P^1$ is big or anti-big exactly when its degree is nonzero. This is the condition $b_1 + b_2 = 2$, and the normal bundle trichotomy follows with $\{b_1, b_2\} = \{1,1\}$, $\{2, 0\}$, or $\{-k, 2+k\}$.
\end{proof}

\begin{remark}\label{rem:wall-geometry}
    The three cases have distinct local geometry. In case (i), the relation $v_1 + v_2 = u_1 + u_2$ shows that $\operatorname{conv}(v_1, v_2)$ and $\operatorname{conv}(u_1, u_2)$ meet at $p = \frac{v_1+v_2}{2}$, so the four rays span the convex cone $\operatorname{cone}(u_1, u_2, v_1, v_2)$, of which $\{\sigma_1, \sigma_2\}$ and $\{\operatorname{cone}(v_1,v_2,u_1), \operatorname{cone}(v_1,v_2,u_2)\}$ are the two simplicial subdivisions without new rays. The wall therefore supports the bistellar move exchanging these triangulations, the toric model of the Atiyah flop \cite{cox2011toric}*{Example 11.1.12 and \S 15.3}. In case (ii), the relation $v_1 + v_2 = 2u_1$, after reordering so that $(b_1, b_2) = (2,0)$, is a linear dependence among $v_1, v_2, u_1$, which span a plane in which $u_1$ generates a ray between those of $v_1$ and $v_2$: within that plane, it is the wall relation of a $(-2)$-curve on a smooth toric surface. In case (iii), with $b_2 < 0$ say, the point $p = \frac{v_1 + v_2}{2} = \frac{b_1 u_1 + b_2 u_2}{2}$ lies in $\operatorname{aff}(T) \setminus T$, so the segment $\operatorname{conv}(v_1, v_2)$ meets the plane $\operatorname{span}(\tau)$ outside $\tau$, and $\sigma_1 \cup \sigma_2$ is not convex.
\end{remark}

\begin{remark}\label{rem:non-invariant-obstruction}
    In dimension at least $3$, an obstructing subvariety need not contain an obstructing invariant curve, and the obstruction can be carried by an invariant subvariety of higher dimension. Let $X$ and $\mathscr{L}$ be as in \Cref{ex:diagonal}. Then we can consider the projective completion of the line bundle $\mathscr{L}^\vee\otimes \omega_X$, which is also toric. The zero section would then have normal bundle as $\mathscr{L}^\vee\otimes \omega_X$, and the canonical bundle of the projective completion restricted to the zero section is precisely given by $\mathscr{L}$, so this is an obstructing surface that doesn't contain any obstructing invariant curve. We do not pursue the classification of obstructing invariant surfaces here.
\end{remark}
\subsection{Obstructing subvarieties and non-standard autoequivalences}
In this section, we work over an algebraically closed field. 
\begin{notation}
     We say an autoequivalence $\Phi$ of $\perf(X)$ or $\D^b(X)$ is \textbf{standard} if $\Phi \cong  f_*(\ecal L \otimes_X -)[n]$ for some $f \in \aut(X)$, $\ecal L \in \pic(X)$ and $n \in \bb Z$. In particular, if $\Phi$ is a standard autoequivalence, then $\Phi$ sends any simple skyscraper sheaf at a closed point to such a sheaf up to shift.  
\end{notation}
When $\omega_X$ is $\otimes$-generating, any autoequivalence of $\perf(X)$ or $\D^b(X)$ is standard by \cite{itoolander2026botamestack}*{Theorem 2.22, Corollary 2.26}, so a natural question is the following:
\begin{question} \label{q:non-standard autoeq}
    If there is an obstructing (invariant) subvariety $Z$ of a toric variety $X$, then is there a non-standard autoequivalence of $\perf(X)$ "associated" to $Z$? 
\end{question}

While the answer is likely to be negative in general, it is positive in the case of smooth projective toric surfaces again by \cite{ito2025derived}*{Lemma 5.7}.
\begin{example}
    Let $X$ be a smooth projective toric surface. Recall that a subvariety is obstructing if and only if it is a $(-2)$-curve $C$. Moreover, note that the spherical twist associated to the structure sheaf $\ecal O_C$ of a $(-2)$-curve is non-standard. See \cite{Uehara2017ATF}*{Example 4.2} for example. 
\end{example}
In higher dimension, an obstructing subvariety is not necessarily invariant, but in this paper, let us focus on invariant ones. For threefolds, we have the following affirmative answer to \Cref{q:non-standard autoeq}.
\begin{example}
    Let $X$ be a smooth complete toric threefold. If there is an obstructing invariant curve $C$ with $\ecal N_{C/X} \cong \ecal O_C(-1) \oplus \ecal O_C(-1)$ in the trichotomy of \Cref{prop:wall-types-smooth}, then the spherical twist associated to $\ecal O_C(-1)$, also known as the flop-flop autoequivalence (cf. \cite{addington2019mukai}*{Theorem A}), provides a non-standard autoequivalence. 
\end{example}
In the rest of the section, we provide more instances of obstructing subvarieties that produce non-standard autoequivalences, using the notion of EZ-spherical twists introduced by Horja \cite{horjatwist}. We work over an algebraically closed field to follow the convention of \cite{horjatwist}.  First, let us recall the following well-known fact.  
\begin{lemma} \label{lemma: Iitaka is derived contraction}
    Let $X = X_\Sigma$ be a complete toric variety and let $\ecal L$ be a nef (equivalently, semiample) line bundle of non-negative Iitaka dimension. Then, the Iitaka fibration $p: X \to B$ is a derived contraction, i.e., $\bb R p_* \ecal O_X \cong \ecal O_B$.
\end{lemma}
\begin{proof}
    By construction (\cite{lazarsfeld2017positivity}*{Theorem 2.1.27, Lemma 2.1.28}, \cite{cox2011toric}*{Theorem 6.2.8}), the Iitaka fibration $p: X \to B$ can be chosen to be a proper surjective toric morphism with connected fibers over a projective toric variety $B$ (see also the proof of \cite{Fujino2007Multiplication}*{Theorem 4.1}). Therefore, $p: X \to B$ is a derived contraction by \cite{Fujino2007Multiplication}*{Lemma 4.3}.  
\end{proof}

\begin{construction} \label{construction: EZ setup for divisors}
    Let $X$ be a smooth complete toric variety and suppose that $E \subset X$ is an obstructing $T$-invariant prime divisor in $X$ (which is necessarily smooth and complete). Suppose that $\omega_X|_E$ or $\omega_X\inv|_E$ is nef and fix $\epsilon \in \{\pm 1\}$ so that $\omega_X^\epsilon|_E$ is nef. Consider the diagram
    \begin{center} \DisableQuotes
        \begin{tikzcd}
         E \arrow[r, hook, "i"] \arrow[d, "p"'] & X \\
         Z                                 &  
    \end{tikzcd}
    \end{center}
    where $p:E \to Z$ is the toric Iitaka fibration with respect to $\omega_X^\epsilon|_E$ as above. By construction, we have 
    \[
    \omega_X^\epsilon|_E = p^* \ecal O_Z(1)
    \]
    for an ample line bundle $\ecal O_Z(1)$ on $Z$, and $\bb Rp_* \ecal O_E \cong \ecal O_Z$ by Lemma \ref{lemma: Iitaka is derived contraction}. 
    Therefore, if $Z$ is smooth and $p$ is flat, then $\ecal O_E$ is an EZ-spherical object with $\theta = \ecal O_Z (-\epsilon)\otimes \omega_Z$ in the sense of Horja (cf. \cite{horjatwist}*{Definition 2.1}), as $\operatorname{codim}_X(E) = 1$ and hence the only requirement is $\bb R p_* \bb R\cal Hom_E(\ecal O_E, \ecal O_E) \cong \ecal O_Z$. In particular, the functor
    \[
    F:=\bb Ri_* \bb L p^*: \perf(Z) \to \perf(X)
    \]
    is a spherical functor (see also \cite{halpern2016autoequivalences}*{Remark 3.7}) and the associated twist $T_F$ gives an autoequivalence of $\perf (X)$ by \cite{horjatwist}*{Theorem 2.11}. In particular, letting $F^R = \bb Rp_* i^!$ denote the right adjoint of $F$, we have the distinguished twisting triangle
    \[
    FF^R(\ecal F) \to \ecal F \to T_F(\ecal F)
    \]
    for any $\ecal F \in \perf X$. Here, since the Grothendieck duality holds in our setup, we do not need the projectivity assumption for $X, E, Z$ (cf. \cite{horjatwist}*{Remark 2.13}), and the other assumptions clearly hold in our situation. Note in particular that we have \[
        i^! \cong \bb R \cal Hom_X(i_*\ecal O_E, -) \cong  \bb L i^* (- \otimes_X \ecal O_X(E) [-1])
        \]
        as $i$ is an inclusion of an effective Cartier divisor (cf. \cite{stacks-project}*{\href{https://stacks.math.columbia.edu/tag/0AU3}{Tag 0AU3 (6), (7)}}).
\end{construction}
\begin{prop}\label{prop-obstructionproduceseztwist}
    Let $X$ be a smooth complete toric variety. If $E \subset X$ is an obstructing invariant prime divisor and if $\omega_X|_E$ or $\omega_X\inv|_E$ is nef and the target $Z$ of the toric Iitaka fibration $p: E \to Z$ has dimension $\leq 1$, then there exists a non-standard autoequivalence $\Phi$ of $\perf X$ such that the following property $(*)$ holds for each closed point $x \in X$:
    \[
    (*) \ \ \text{if $x \not \in E$, then $\Phi(\kappa(x)) = \kappa(x)$ and if $x \in E$, then $\dim \supp \Phi(\kappa(x)) > 0$.} \qedhere 
    \]
\end{prop}
  We treat the cases $\dim Z = 0$ and $\dim Z = 1$ separately, whence the proof of \Cref{prop-obstructionproduceseztwist} follows.

    \begin{lemma}
        If $\dim Z = 0$, then $\omega_X|_E \cong \ecal O_E$ and the corresponding twist agrees with the spherical twist with respect to the spherical object $i_*\ecal O_E$. In particular, the corresponding spherical twist is non-standard and satisfies the property $(*)$.
    \end{lemma}
    \begin{proof}
        When $\dim Z = 0$, we see $Z \cong \spec k$ and $p$ is flat. Moreover, $\omega_X|_E = p^* \ecal O_{Z} = \ecal O_E
        $ by construction. In particular, $i_*\ecal O_E$ is a spherical object by \cite{horjatwist}*{Remark 2.3}. Recall that the last claim follows for example by \cite{Uehara2017ATF}*{Example 4.2 (ii)}, noting that the arguments there follow for spherical twists on any variety. Alternatively, apply the same proof as in the next lemma. 
    \end{proof}
    \begin{lemma}
        If $\dim Z = 1$, then $Z \cong \bb P^1$ and $p: E \to Z$ is flat. Moreover, the corresponding EZ twist is non-standard and satisfies the property $(*)$. 
    \end{lemma}
    \begin{proof}
        Since $Z$ needs to be a projective toric curve, $Z \cong \bb P^1$ and $p:E \to Z$ is flat for example by \cite{Har77}*{III.9.7}. Therefore, it suffices to show the corresponding EZ twist satisfies the property $(*)$. As in \Cref{construction: EZ setup for divisors}, 
    let $F= \bb Ri_* \bb Lp^*$ and its right adjoint $F^R = \bb Rp_* i^!$. First, take $x \in X \setminus E$. Then, since $i^!\kappa(x) = 0$, we have $T_F(\kappa(x)) \cong \kappa(x)$ by the twisting triangle. Next, take $x \in E$. Since $i^! \kappa(x) \cong \bb Li^*\kappa(x)[-1]$ is a non-zero perfect complex supported at $x \in E$, we see that $F^R(\kappa(x))$ is at most supported at $p(x)$. It is indeed non-zero since
    \[
    \hom_Z(\ecal O_Z, F^R\kappa(x)) = \hom_X(F(\ecal O_Z), \kappa(x)) = \hom_X(i_*\ecal O_E, \kappa(x)) = \hom_E(\ecal O_E, \kappa(x)) \neq 0.
    \]
    Therefore, since $p$ is flat and $i$ is a closed immersion, we can compute $\supp FF^R\kappa(x)$ by applying $i_*p^*$ to the cohomology sheaves and therefore $\supp FF^R\kappa(x) = p\inv(p(x))$. Since the support is closed and  $\dim p\inv (p(x)) = \dim E - 1 > 0$, the twist triangle shows that $\supp FF^R\kappa(x) = p\inv(p(x)) \subset \supp T_F\kappa(x) \cup \{x\}$ and hence $\dim \supp T_F \kappa(x) \geq \dim p\inv (p(x))> 0$. Hence, we are done.
    \end{proof}

 When $\dim Z = 2$, the Iitaka fibration $p: E \to Z$ is not necessarily flat and the target $Z$ is not necessarily smooth, so we cannot directly apply Horja's construction to construct a twist autoequivalence. However, it would be interesting to see exactly when an obstructing subvariety to which the canonical bundle restricts as a nef or anti-nef line bundle gives rise to a twist-like autoequivalence. Moreover, it would also be interesting to see when obstructing subvarieties correspond to a suitable wall in VGIT whose associated window shift is a non-standard autoequivalence in the setting of \cite{halpern2016autoequivalences}. Another natural question is whether the data of the obstructing subvariety under the toric mirror symmetry provides a way to construct a non-standard autoequivalence on the other side. We leave these directions as future research.   

\bibliography{bib}
\end{document}